\documentclass[12pt]{article}
\usepackage[T1]{fontenc}
\usepackage{lmodern}
\usepackage[margin=1in]{geometry}
\usepackage{amsmath,amssymb,amsthm,mathtools}
\usepackage{microtype,needspace}
\usepackage{titlesec}
\titleformat{\section}[block]
  {\normalfont\normalsize\scshape\centering}
  {\thesection.}{0.5em}{}
\titlespacing*{\section}{0pt}
  {0.7\baselineskip plus \baselineskip}{0.5\baselineskip}
\usepackage{tikz,pgfplots}
\usetikzlibrary{arrows.meta,calc,decorations.pathreplacing}
\pgfplotsset{compat=1.18}
\usepackage{booktabs,tabularx,array,caption,float}
\usepackage[section]{placeins}
\pgfplotsset{paperplot/.style={
 axis lines=left,tick align=outside,
 axis line style={line width=.4pt},tick style={line width=.4pt},
 ticklabel style={font=\small},label style={font=\small},
 every axis plot/.append style={line width=.9pt},scaled ticks=false}}
\tikzset{
 indexvertex/.style={circle,draw,fill=white,inner sep=0pt,
 minimum size=5.3mm,font=\small,line width=.5pt},
 matrixedge/.style={-{Stealth[length=1.5mm]},line width=.8pt,
 shorten <=2.8mm,shorten >=2.8mm},
 pairingedge/.style={densely dashed,line width=.8pt},
 cycleedge/.style={-{Stealth[length=1.6mm]},line width=.7pt}}
\usepackage[hidelinks,hypertexnames=false]{hyperref}
\allowdisplaybreaks[0]
\newtheorem{theorem}{Theorem}[section]
\newtheorem{lemma}[theorem]{Lemma}
\newtheorem{proposition}[theorem]{Proposition}
\newtheorem{corollary}[theorem]{Corollary}
\newtheorem{mainthm}{Theorem}

\theoremstyle{definition}

\newtheorem{example}[theorem]{Example}

\theoremstyle{remark}
\newtheorem{remark}[theorem]{Remark}

\DeclareMathOperator{\Tr}{Tr}
\DeclareMathOperator{\Var}{Var}
\DeclareMathOperator{\Cov}{Cov}
\DeclareMathOperator{\Pois}{Pois}
\DeclareMathOperator{\sgn}{sgn}

\newcommand{\R}{\mathbb R}

\newcommand{\op}{\mathrm{op}}
\newcommand{\F}{\mathrm F}
\newcommand{\dto}{\Rightarrow}
\newcommand{\Pp}{\mathbb P}
\newcommand{\Ee}{\mathbb E}
\newcommand{\cP}{\mathcal P}
\newcommand{\indep}{\mathrel{\perp\!\!\!\perp}}

\title{Conditional Gaussian Limits and Character Laws\\for Randomized Symmetric-Group Representations}
\hypersetup{pdftitle={Conditional Gaussian Limits and Character Laws for Randomized Symmetric-Group Representations},
 pdfauthor={Jonathan Guo},pdfsubject={Explicit character laws and Haar-conjugation limits},
 pdfkeywords={Young diagram, symmetric group, character polynomial, Haar orthogonal, randomized basis}}
\date{}
\begin{document}
\author{Jonathan Guo}
\maketitle
\begin{abstract}
Let $\rho_n$ be a real orthogonal representation of the symmetric group
$S_n$ of dimension $d_n\to\infty$, let $\sigma_n$ be a uniform random
permutation, and let $O_n$ be an independent Haar-orthogonal matrix. We
study the linear statistics $\Tr(A_nO_n\rho_n(\sigma_n)O_n^T)$ for
deterministic real test matrices $A_n$ normalized only in Frobenius
norm, which allows tests that isolate individual entries. Conditionally
on $\sigma_n$, the statistic has mean $s_n\chi_n(\sigma_n)$, where
$\chi_n$ is the character of $\rho_n$ and $s_n=d_n^{-1}\Tr A_n\to s$.
We show that the centered statistic has conditional Gaussian limit
$N(0,1-s^2)$ for every such test if and only if the trivial and sign
constituents of $\rho_n$ occupy $o(d_n)$ dimensions, and that this
limit is asymptotically independent of every convergent
$\sigma_n$-measurable observable. The key input is a tensor-square
estimate, based on the Frobenius--Schur identity and the minimum degree
of non-linear irreducible representations of $S_n$, which shows that
this single condition controls both normalized traces needed in the
Haar-conjugation argument. For the uncentered statistic, convergence in
distribution of $\chi_n(\sigma_n)$ to a law $\Xi$ is necessary and
sufficient for every test statistic to converge to
$s\Xi+\sqrt{1-s^2}\,G$, with $G$ standard Gaussian and independent of
$\Xi$. For every irreducible sequence with $d_n\to\infty$, finitely
many entries scaled by $\sqrt{d_n}$ have independent standard Gaussian
limits. Finally, we compute $\Xi$ explicitly for diagrams with a fixed
first-row tail, their transposes, and diagrams with sublinear largest
hook, answering in these cases a question of Tsou.

\vspace{12pt}
\end{abstract}

\section{Introduction}\label{sec:introduction}

A change of basis can change the entries of a representation matrix, but
it does not change its trace. Thus we should not expect all information
about the representation to disappear when the basis is randomized.
Even if the centered fluctuation becomes Gaussian, the trace can still
contribute to the limit. We want to understand this contribution and,
when possible, read it from the Young diagrams rather than assume a
limit for the random traces.

For a uniform permutation matrix $P_n$ and an independent Haar-orthogonal
matrix $O_n$, Tsou \cite[Theorem~3.2]{TsouRandomized} proves that
$\Tr(A_nO_nP_nO_n^T)\dto sY+\sqrt{1-s^2}\,G$ for deterministic real
tests with $\Tr(A_nA_n^T)=n$ and $n^{-1}\Tr A_n\to s$. Here
$Y\sim\Pois(1)$ and $G\sim N(0,1)$ are independent. His Remark~3.2
asks whether similar distributional limits hold for linear combinations
of entries of higher-dimensional representations of $S_n$. We approach
this question in two parts. First we prove a conditional Gaussian
limit for the centered fluctuation. Then we identify the character
that appears in the conditional mean. Theorems~\ref{thm:A}--\ref{thm:C}
below, stated after the necessary notation, summarize the answer.

Suppose $\rho_n:S_n\to\mathrm O(d_n)$ is a real orthogonal
representation. Let $\sigma_n$ be uniform on $S_n$, and independently
let $O_n$ be Haar on $\mathrm O(d_n)$. We use the notation
\[
 Q_n=\rho_n(\sigma_n),\qquad
 \chi_n=\Tr\rho_n,\qquad
 C_j^{(n)}=C_j(\sigma_n),\qquad \varepsilon_n=\sgn(\sigma_n).
\]
Here $n$ is the degree of the symmetric group, whereas $d_n$ is the
dimension of the representation. These are different quantities: both
tend to infinity in the representation results, but they need not grow
at the same rate. We write $\|A\|_\F=(\Tr AA^T)^{1/2}$ for the
Frobenius norm and $\|A\|_\op$ for the Euclidean operator norm.

Unless stated otherwise, our test is a deterministic real
$d_n\times d_n$ matrix $A_n$ satisfying
\begin{equation}
 \Tr(A_nA_n^T)=d_n(1+o(1)),\qquad
 s_n=\frac{\Tr A_n}{d_n}\longrightarrow s.
 \tag{A}\label{eq:test-normalization}
\end{equation}
By Cauchy--Schwarz, $|s|\le1$. Notice that this normalization allows
both $A_n=I_{d_n}$ and $A_n=\sqrt{d_n}\,E_{ba}$, where $E_{ba}$ is a
matrix unit. The latter has operator norm $\sqrt{d_n}$, so a uniform
operator-norm bound would exclude the tests used to study individual
entries.

Let $T_n=\Tr(A_nO_nQ_nO_n^T)$. Once we condition on $\sigma_n$, the
mean over the random basis is determined by the character. We can
therefore write
\begin{equation}
 T_n=s_n\chi_n(\sigma_n)+X_n,
 \qquad \Ee[X_n\mid\sigma_n]=0.
 \label{eq:intro-decomposition}
\end{equation}
The first term does not depend on the choice of basis. The second is
what remains after subtracting the conditional mean.

Let $t_n$ denote the total multiplicity of the trivial and sign
representations in $\rho_n$. Both are one-dimensional: their Young
diagrams are a single row and a single column, respectively. Our first
result characterizes the centered Gaussian limit in terms of $t_n$
alone.

\begin{mainthm}[Theorem~\ref{thm:universality}]\label{thm:A}
Let $\rho_n:S_n\to\mathrm O(d_n)$ be arbitrary real representations
with $d_n\to\infty$. For a deterministic test satisfying
\eqref{eq:test-normalization}, let $X_n=T_n-s_n\chi_n(\sigma_n)$ be the
centered statistic in \eqref{eq:intro-decomposition}. Then
\[
 \frac{t_n}{d_n}\longrightarrow0
 \quad\Longleftrightarrow\quad
 X_n\dto N(0,1-s^2)\ \text{for every such test sequence.}
\]
Under this condition, moreover,
$\Ee[e^{iuX_n}\mid\sigma_n]\to\exp(-u^2(1-s^2)/2)$ in probability and
in $L^1$ for every $u\in\R$, and the Gaussian fluctuation is jointly
asymptotically independent of any convergent fixed-dimensional
$\sigma_n$-measurable vector.
\end{mainthm}

The condition is necessary even if we require only unconditional
convergence. There is no restriction on the shapes or multiplicities
of the other constituents, and the unnormalized character need not
converge or even be tight. There is also an entry conclusion that does
not require a character limit (Corollary~\ref{cor:all-irreducible-entries}):
for every irreducible sequence $S^{\lambda_n}$ with
$f^{\lambda_n}\to\infty$, any fixed collection of distinct entries of
$O_n\rho_n(\sigma_n)O_n^T$, scaled by $\sqrt{d_n}$, has jointly
independent standard Gaussian limits, with no stabilization or
geometric restriction on $\lambda_n$.

The uncentered statistic carries the character through the first term
in \eqref{eq:intro-decomposition}. Our second result says that
convergence of this character is the only additional thing to check.

\begin{mainthm}[Theorem~\ref{thm:character-equivalence}]\label{thm:B}
Let $\rho_n:S_n\to\mathrm O(d_n)$ be arbitrary real representations
with $d_n\to\infty$; no multiplicity hypothesis is imposed. For a real
random variable $\Xi$ with a fixed law,
$\chi_n(\sigma_n)\dto\Xi$ if and only if, for every deterministic test
satisfying \eqref{eq:test-normalization},
\[
 T_n\dto s\Xi+\sqrt{1-s^2}\,G,\qquad G\sim N(0,1),\quad G\indep\Xi.
\]
\end{mainthm}

What remains is to identify $\Xi$ from the Young diagrams. The third
result does this for three families of diagrams and their direct sums.
For a fixed partition $\mu$, write $\mu[n]$ for the partition of $n$
obtained by placing a first row of length $n-|\mu|$ above $\mu$, and
write $\lambda'$ for the transpose of $\lambda$.

\begin{mainthm}[Theorem~\ref{thm:diagram} and Corollary~\ref{cor:irreducible}]
\label{thm:C}
Suppose $d_n\to\infty$ and that, for all sufficiently large $n$,
\[
 \rho_n\cong
 \bigoplus_{\mu\in\mathcal T}(S^{\mu[n]})^{\oplus a_\mu}
 \ \oplus\!
 \bigoplus_{\mu\in\mathcal T}(S^{(\mu[n])'})^{\oplus b_\mu}
 \ \oplus R_n,
\]
where $\mathcal T$ is a finite set of fixed partitions, $a_\mu,b_\mu$
are fixed nonnegative integers, and $R_n$ is any real representation
of $S_n$ whose constituent diagrams have largest hook $o(n)$, with no
restriction on its multiplicities. Let $P_1,P_2,\ldots,\varepsilon,G$
be independent, with $P_j\sim\Pois(1/j)$,
$\Pp(\varepsilon=1)=\Pp(\varepsilon=-1)=1/2$, and $G\sim N(0,1)$, and
let $q_\mu$ be the explicit polynomial in the cycle counts
$x_1,\ldots,x_{|\mu|}$ defined in \eqref{eq:stable-polynomial}. Then
\[
 \chi_n(\sigma_n)\dto
 \Xi:=\sum_{\mu\in\mathcal T}a_\mu q_\mu(P)
      +\varepsilon\sum_{\mu\in\mathcal T}b_\mu q_\mu(P),
\]
and consequently $T_n\dto s\Xi+\sqrt{1-s^2}\,G$ for every test
satisfying \eqref{eq:test-normalization}. In particular, for a single
irreducible representation $S^{\lambda_n}$, the limit $\Xi$ equals
$q_\mu(P)$ when $\lambda_n=\mu[n]$, equals $\varepsilon q_\mu(P)$ when
$\lambda_n=(\mu[n])'$, and equals $0$ when
$\lambda_{n,1}+\lambda'_{n,1}=o(n)$.
\end{mainthm}

All summands in $\Xi$ use the same Poisson variables and the same sign,
because they come from one permutation. An arbitrary small-hook
remainder does not change the limiting character. These diagram
conditions describe the non-Gaussian part only; they are not needed
for the centered Gaussian theorem. Taking $\rho_n$ to be the defining
permutation representation $S^{(n)}\oplus S^{(n-1,1)}$ gives
$\Xi=P_1$ and recovers Tsou's Poisson-plus-Gaussian law. Thus
Theorem~\ref{thm:C} answers the question in Tsou's Remark~3.2 for
fixed first-row tails, their transposes, and small-hook diagrams,
while Theorem~\ref{thm:A} answers the centered form of the question
for arbitrary real representations.

There is a stronger conclusion when the character is given by one
fixed cycle polynomial at every sufficiently large degree, rather
than merely converging in distribution.
Proposition~\ref{prop:exact-tail-equivalence} uses
Murnaghan--Nakayama to prove the converse: with a fixed optional sign
factor, such a formula forces the diagram to stabilize below the
first row or outside the first column.
Corollary~\ref{cor:row-column-polynomials} gives the same alternatives
when both polynomial terms are allowed in an irreducible character.

We now describe the arguments. For arbitrary random orthogonal matrices $Q_n$,
Theorem~\ref{thm:linear} gives a conditional Gaussian limit of
variance $1-s^2$ when the two normalized traces $\Tr Q_n/d_n$ and
$\Tr(Q_n^2)/d_n$ tend to zero in probability. The analytic input is a
comparison of each fixed Haar moment of $\Tr(AOBO^T)$ with the
corresponding moment of a Gaussian quadratic form, uniformly over all
Frobenius-normalized tests (Lemma~\ref{lem:weingarten-moment-comparison}
and Theorem~\ref{thm:weingarten-cumulants}). When the matrices come
from representations of $S_n$, we replace these two hypotheses by one
condition on the irreducible constituents. The estimate that relates
the two normalized characters is
\[
 \Ee\!\left[\left(\frac{\chi_n(\sigma_n^2)}{d_n}\right)^2\right]
 \le 2\Ee\!\left[\left(\frac{\chi_n(\sigma_n)}{d_n}\right)^2\right]
       +\frac1{n-1}\qquad(n\ge5).
\]
We prove this in Proposition~\ref{prop:square-character-bound} by
decomposing the tensor square and using the Frobenius--Schur identity
together with the minimum-degree bound for non-linear irreducibles
\cite{Rasala,BNOT}. This gives the sufficiency in
Theorem~\ref{thm:A}. To prove necessity at the level of distributions,
we use a skew-symmetric test. A uniform fourth-moment bound then allows
us to pass from convergence in distribution to convergence of second
moments.

For Theorem~\ref{thm:B}, tightness of the ordinary character, together
with the estimate above, makes both normalized traces tend to zero.
Thus there is no separate condition to check at the square.
Conversely, taking $A_n=I_{d_n}$ returns the character itself.

For Theorem~\ref{thm:C}, the fixed-tail character formula is classical
\cite[Example~I.7.14(b)]{Macdonald}; see also Garsia--Goupil
\cite{GarsiaGoupil} and Narayanan--Paul--Prasad--Srivastava
\cite{NPPS}. For broader background on asymptotic characters, see
Vershik--Kerov \cite{VershikKerov} and the partial-permutation algebra
of Ivanov--Kerov \cite{IvanovKerov}. Once we have a fixed character
polynomial, its limiting law follows by substituting the Poisson
limits of the small cycle counts. The small-hook case uses a different
observation: one long-cycle event makes every constituent character
zero at both the permutation and its square.

We next place these results relative to existing work. For the
Gaussian part, we use Haar-conjugation fluctuation theory.
Related unitary statistics $\Tr(AUBU^*)$ occur in spherical-integral
(Harish-Chandra--Itzykson--Zuber) asymptotics; see Guionnet--Ma\"ida
\cite{GuionnetMaida} and the Haar integration formulas of
Collins--\'{S}niady \cite{CollinsSniady}. Relevant second-order-freeness
results include Mingo--\'{S}niady--Speicher
\cite[Theorem~3.15]{MingoSniadySpeicher} and Mingo--Popa
\cite[Proposition~52]{MingoPopa}; the real theory also includes
transpose terms \cite{RedelmeierReal}. The unitary results of
Benaych-Georges--C\'ebron--Rochet
\cite[pp.~5--6, Hypothesis~1(d), Theorems~2.1 and~2.3, and Remark~2.5]{BCR}
are closest to the setting here. Their normalized pair-trace conditions
allow entry-scale tests. Their trace-free Theorem~2.1 does not require
a limit law for the unnormalized traces. Their Theorem~2.3 assumes
convergence of the centered traces through Hypothesis~2, and gives a
trace component independent of the Gaussian fluctuation. They also
consider uniform permutation matrices under unitary conjugation
\cite[p.~8, Corollary~2.9]{BCR}.

The present paper differs from these results in four respects.
\begin{enumerate}
\item[(i)] \emph{Orthogonal rather than unitary conjugation.}
Orthogonal conjugation need not produce unitary invariance, and the
real moment expansion contains additional transpose contractions
(Section~\ref{sec:transfer}). The orthogonal laws therefore do not
follow directly from the unitary ones, even for fixed tails.
\item[(ii)] \emph{Frobenius-only normalization of the test.}
We must allow $A_N=\sqrt N E_{ba}$. Its operator norm is $\sqrt N$, so
it does not satisfy a bounded-operator-norm hypothesis. It also fails
the polynomial trace-moment and trace-cumulant hypotheses of the cited
second-order formulations: already $N^{-1}\Tr((A_NA_N^T)^2)=N$. This
is why we work with estimates that require only Frobenius
normalization.
\item[(iii)] \emph{A conditional statement without tightness.}
We prove a conditional orthogonal limit, with independence from every
convergent observable of the underlying randomness, without requiring
tightness or convergence of the unnormalized trace.
\item[(iv)] \emph{A multiplicity criterion for $S_n$.}
For symmetric-group representations, the two normalized-trace
conditions of the Haar argument are replaced by the single condition
$t_n=o(d_n)$, which is also necessary. Identifying the remaining
character law is then a separate question about the diagrams.
\end{enumerate}
The deterministic Theorem~\ref{thm:deterministic-clt} is proved by a
self-contained Weingarten comparison whose coefficient estimate is a
coarse fixed-order consequence of Collins--Matsumoto
\cite[Theorem~4.11]{CollinsMatsumoto}. What we need from its proof is
the uniformity over all Frobenius-normalized tests, which the
conditional argument of Theorem~\ref{thm:linear} requires.

Tsou's theorem already gives an independent Poisson and Gaussian
component \cite[Theorem~3.2]{TsouRandomized}. Following Tsou, we also
acknowledge P.~Dey's moment-method approach, recorded there as a 2008
personal communication. The centered criterion here applies to
arbitrary real representations, and the entry result applies to every
growing-dimensional irreducible sequence. In the families considered
below, the explicit character calculations then give the uncentered
limits.

We begin in Section~\ref{sec:transfer} with the Haar argument. We then
prove the square-character estimate and use it to obtain the
multiplicity criterion, the equivalence with character convergence,
and the entry limits. In Section~\ref{sec:characters}, we derive the
fixed-tail polynomials, prove the exact row-or-column converse, and
compute limiting character laws using Poisson cycle counts and the
small-hook argument. We finish with examples and discuss which
character laws are not identified by these calculations. Each proof
appears with its statement.

\newpage

\section{Gaussian limits and the multiplicity criterion}\label{sec:transfer}
\label{sec:haar-proofs}\label{sec:notation}\label{sec:weingarten}

We first set aside the representation and consider an arbitrary
orthogonal matrix. Starting with a deterministic matrix lets us see
exactly where the trace conditions enter. We compute the centered
variance, compare the Haar moments with Gaussian moments, and then
allow the matrix to be random by conditioning. The character
calculations will come afterward.

To see why two normalized traces appear, first consider a single
entry. For deterministic $B\in\mathrm O(N)$, let
$\beta=N^{-1}\Tr B$, $\tau=N^{-1}\Tr(B^2)$, and $D=B-\beta I_N$.
We subtract the scalar part before making a Gaussian approximation,
so that the exact orthogonality of the rows of $O$ is respected.
Replacing a fixed number of rows by independent vectors
$g_a\sim N(0,I_N/N)$ suggests the centered Gaussian expression
$Z_{ab}=\sqrt N\,g_a^TDg_b$. Gaussian pairings give
\begin{align*}
 \Var Z_{ab}&=1-\beta^2+\mathbf1_{\{a=b\}}(\tau-\beta^2),\\
 \Cov(Z_{ab},Z_{ba})&=\tau-\beta^2\qquad(a\ne b).
\end{align*}
Notice that the second trace appears both in the diagonal variance
and in the covariance of transposed entries. This explains its role
at the level of individual entries; Remark~\ref{rem:trace-sharpness}
below shows that both traces are needed for the variance limit to
hold uniformly over tests. To handle an arbitrary normalized
test, we must control the calculation uniformly over all its
coefficients. The Weingarten contraction bounds below will do this.

For notation, write $\operatorname{tr}_N=N^{-1}\Tr$ and
$\langle A,B\rangle_\F=\Tr(AB^T)$. It will be useful to separate a
matrix into its three Frobenius-orthogonal components. Let
$A_{s,0}=(A+A^T)/2-\operatorname{tr}_N(A)I_N$ and
$A_a=(A-A^T)/2$. Then
$A=\operatorname{tr}_N(A)I_N+A_{s,0}+A_a$: the first part is scalar,
the second is symmetric and traceless, and the third is skew-symmetric.
With this notation, the standard Haar mean and variance identities
take the following form \cite{CollinsSniady,GoodmanWallach,Meckes}.

\begin{proposition}\label{prop:haar-mean}
\label{prop:exact-variance}\label{lem:orthogonal-parts}
For real $N\times N$ matrices $A,B$, $N\ge2$, and Haar $O\in\mathrm O(N)$,
\[
 \Ee_O[OBO^T]=\operatorname{tr}_N(B)I_N,\qquad
 \Ee_O\Tr(AOBO^T)=\operatorname{tr}_N(A)\Tr B,
\]
and
\begin{equation}
 \Var_O\Tr(AOBO^T)=
 \frac{\|A_{s,0}\|_\F^2\|B_{s,0}\|_\F^2}{d_s}
 +\frac{\|A_a\|_\F^2\|B_a\|_\F^2}{d_a},
 \label{eq:exact-variance}
\end{equation}
where $d_s=(N-1)(N+2)/2$ and $d_a=N(N-1)/2$.
When $B$ is orthogonal,
\begin{equation}
 \|B_{s,0}\|_\F^2=\frac{N+\Tr(B^2)}2-\frac{(\Tr B)^2}{N},
 \qquad \|B_a\|_\F^2=\frac{N-\Tr(B^2)}2.
 \label{eq:Bs0-norm}
\end{equation}
\end{proposition}

\begin{proof}
We first compute the mean. By Haar invariance, conjugating
$\Ee_O[OBO^T]$ by any orthogonal matrix leaves it unchanged. Thus this
mean is a scalar matrix. Its trace is $\Tr B$, so the scalar must be
$\operatorname{tr}_N(B)$. Taking the trace after multiplication by $A$
gives the second identity for the mean.

For the variance, recall the decomposition into scalar,
symmetric-traceless, and skew-symmetric parts. The scalar part does not
fluctuate under conjugation. The other two parts belong to
nonisomorphic irreducible real orthogonal modules of dimensions $d_s$
and $d_a$. These are the degree-two harmonic symmetric and exterior
tensor representations of the full orthogonal group
\cite[Section~10.2]{GoodmanWallach}. Since the modules are
nonisomorphic, their cross covariance is zero.

Now consider either one of these modules. The averaged rank-one
covariance operator is self-adjoint and commutes with conjugation.
Consequently each of its eigenspaces is invariant. Irreducibility then
means that there can be only one eigenspace, so the operator is a scalar
multiple of the identity. Its trace is $\|B_\bullet\|_\F^2$, and the
module has dimension $d_\bullet$. Hence the scalar is
$\|B_\bullet\|_\F^2/d_\bullet$. Applying this operator to the
corresponding component of $A$, and adding the two contributions, gives
\eqref{eq:exact-variance}.

Finally, suppose $B$ is orthogonal. Expand the squared Frobenius norms
of $(B+B^T)/2$ and $(B-B^T)/2$, using $BB^T=I_N$. For the symmetric
part we must also subtract the squared norm of the scalar component,
which is $(\Tr B)^2/N$. This gives the two formulas in
\eqref{eq:Bs0-norm}.
\end{proof}

\begin{proposition}
\label{prop:variance-error}
For $N\ge2$, $\|A\|_\F^2=N$, $B\in\mathrm O(N)$, and
$a=N^{-1}\Tr A$, $b=N^{-1}\Tr B$, $\tau=N^{-1}\Tr(B^2)$,
\begin{equation}
 \left|\Var_O\Tr(AOBO^T)-(1-a^2)\right|
 \le \frac2N+2|\tau|+2b^2.
 \label{eq:variance-error}
\end{equation}
\end{proposition}

\begin{proof}
Let $x=\|A_{s,0}\|_\F^2/N$ and $y=\|A_a\|_\F^2/N$. Since the
three components of $A$ are orthogonal and $\|A\|_\F^2=N$, we have
$x+y=1-a^2\le1$. Substituting into the exact variance formula gives
\[
 x\frac{N^2(1+\tau-2b^2)}{(N-1)(N+2)}
   +y\frac{N(1-\tau)}{N-1}.
\]
We want to compare this expression with $x+y$. Notice that
$N^2/((N-1)(N+2))$ is at most one and differs from one by at most
$1/N$. Similarly, $N/(N-1)$ is at most two and differs from one by at
most $2/N$. The terms involving $\tau$ therefore contribute at most
$2|\tau|$, and the term involving $b^2$ contributes at most $2b^2$.
The remaining error is at most $2/N$, since $x,y\ge0$ and $x+y\le1$.
Adding these bounds proves the result.
\end{proof}

\begin{corollary}\label{cor:variance-limit}
If $\|A_N\|_\F^2=N$, $a_N=\operatorname{tr}_N A_N\to s$, and
$B_N\in\mathrm O(N)$ satisfies
$\operatorname{tr}_N B_N\to0$ and $\operatorname{tr}_N(B_N^2)\to0$, then
$\Var_O\Tr(A_NO_NB_NO_N^T)\to1-s^2$.
\end{corollary}
\begin{proof}
Apply \eqref{eq:variance-error} with $A=A_N$ and $B=B_N$. Each term on
the right tends to zero by the hypotheses. Thus the variance differs
from $1-a_N^2$ by a quantity tending to zero. Since
$1-a_N^2\to1-s^2$, the desired limit follows.
\end{proof}

\begin{remark}\label{rem:trace-sharpness}
Suppose $B_N\in\mathrm O(N)$ is deterministic, and write
$b_N=\operatorname{tr}_N B_N$ and $\tau_N=\operatorname{tr}_N(B_N^2)$.
The two trace conditions are not just sufficient if we require the
variance limit uniformly over all normalized tests. In fact,
$b_N\to0$ and $\tau_N\to0$ if and only if
\begin{equation}
 \sup_{\|A\|_\F^2=N}
 \left|\Var_O\Tr(AO_NB_NO_N^T)
       -\bigl(1-(\operatorname{tr}_N A)^2\bigr)\right|
 \longrightarrow0.
 \label{eq:uniform-variance-sharpness}
\end{equation}
One direction follows from \eqref{eq:variance-error}. For the other,
consider the two traceless normalized tests
$A_N^-=\sqrt{N/2}(E_{12}-E_{21})$ and
$A_N^+=\sqrt{N/2}(E_{12}+E_{21})$.
By \eqref{eq:exact-variance}--\eqref{eq:Bs0-norm}, their variances are
\[
 \Var_O\Tr(A_N^-O_NB_NO_N^T)=\frac{N(1-\tau_N)}{N-1},\qquad
 \Var_O\Tr(A_N^+O_NB_NO_N^T)
   =\frac{N^2(1+\tau_N-2b_N^2)}{(N-1)(N+2)}.
\]
If the first variance tends to one, then $\tau_N\to0$. Substituting
this into the second variance shows that $b_N\to0$ as well.
The calculation holds for each realization of $B_N$. Thus, when
$B_N$ is replaced by $Q_n$ and $N$ by $d_n$, the same equivalence
holds in probability for conditional variances.

A single test need not detect both traces, even when its
symmetric-traceless part is nonzero. For example, take even $N$ and let
$B_N=\operatorname{diag}(I_{N/2},-I_{N/2})$ and
$A_N=\sqrt N E_{12}$. Then $b_N=0$ and $\tau_N=1$, but
\[
 \Var_O\Tr(A_NO_NB_NO_N^T)
   =\frac{N^2}{(N-1)(N+2)}\longrightarrow1.
\]
Thus the off-diagonal entry has the desired limiting variance even
though the second trace condition fails. Notice, however, that
$O_NB_NO_N^T$ is symmetric. Its transposed entries are equal, so their
joint law still detects the failure of that condition.
\end{remark}

The variance suggests the Gaussian law we should expect, but matching
a variance does not prove convergence to that law. We next compare
each fixed Haar moment with the corresponding moment of a Gaussian
quadratic form. The comparison must be uniform over all
Frobenius-normalized tests, including the scaled matrix units used
for individual entries.

Throughout the next estimates, constants may depend on the fixed
order and on $M$, but not on the matrices or the dimension. The
estimates hold for all sufficiently large $N$, with the threshold
depending only on those parameters. Recall the moment--cumulant
identity
\begin{equation}
 \Ee X^r=\sum_{\pi\in\cP(r)}\prod_{V\in\pi}\kappa_{|V|}(X).
 \label{eq:moment-cumulant}
\end{equation}
Here $\cP(r)$ is the set of set partitions, and $\kappa_j$ denotes
the ordinary $j$th cumulant. If the first cumulant tends to $0$, the
second cumulant tends to $v$, and every fixed higher cumulant tends
to zero, then the moments tend to those of $N(0,v)$. Those moments
determine the Gaussian law, so $X_N\dto N(0,v)$. In the case $v=0$,
convergence of the second moment already suffices.

Let $\mathcal P_2(2k)$ denote the pair partitions and let
$L_k=(2k-1)!!$ be their number. For a pairing $p$, define
$\delta_p(i)=\prod_{\{a,b\}\in p}\mathbf1_{\{i_a=i_b\}}$.
Thus $\delta_p(i)$ imposes the index equalities prescribed by $p$.
For two pairings $p,q$, let $c(p,q)$ count the components of their
union, retaining doubled edges. We have $c(p,p)=k$, since every
component is a doubled edge. If $p\ne q$, at least two pairs belong
to one component, so $c(p,q)\le k-1$.

\begin{lemma}
\label{lem:weingarten-fixed-order}
For fixed $k$ and sufficiently large $N$,
\begin{align}
 \Ee\prod_{a=1}^{2k}O_{i_a j_a}
 &=\sum_{p,q}\delta_p(i)\delta_q(j)\operatorname{Wg}_{N,k}(p,q),
 \label{eq:weingarten-integration}\\
 \left|\operatorname{Wg}_{N,k}(p,q)-N^{-k}\mathbf1_{\{p=q\}}\right|
 &\le w_kN^{-k-1}.
 \label{eq:weingarten-bound}
\end{align}
The coefficient matrix is the inverse of $(N^{c(p,q)})_{p,q}$.
\end{lemma}
\begin{proof}
Recall that averaging $O^{\otimes2k}$ over the orthogonal group gives
the projection onto the invariant tensors. By the orthogonal first
fundamental theorem, the pairing vectors span this invariant space for
every $N$ \cite[Theorem~5.3.3]{GoodmanWallach}. We therefore need to
understand their Gram matrix.

The entry indexed by $p,q$ is $N^{c(p,q)}$. On the diagonal this is
$N^k$, while off the diagonal it is at most $N^{k-1}$. Thus we can write
the Gram matrix as $N^k(I+E)$, where
$\|E\|_\infty\le(L_k-1)/N$ in the row-sum norm. If
$N>2(L_k-1)$, then $\|E\|_\infty<1/2$. Expanding the inverse of
$I+E$ as a geometric series shows that the inverse Gram matrix differs
from $N^{-k}I$ by at most $2(L_k-1)N^{-k-1}$ in this norm. The formula
for projection in terms of a Gram matrix now gives both assertions.

The lower bound on $N$ is only used to obtain this inverse estimate.
The pairing vectors still span the invariant space for smaller $N$,
although they may then be linearly dependent. When $k=1$, the formula
is exact. See also \cite{CollinsSniady} and the sharper estimates in
\cite[Theorem~4.11]{CollinsMatsumoto}.
\end{proof}

\begin{lemma}
\label{lem:weingarten-moment-comparison}
Suppose
\begin{equation}
 \Tr C=\Tr D=0,\qquad \|C\|_\F\le\sqrt N,\qquad\|D\|_\op\le M.
 \label{eq:weingarten-CD}
\end{equation}
For Haar $O$ and a matrix $G$ with independent $N(0,1/N)$ entries, set
$Y=\Tr(CODO^T)$ and $V=\Tr(CGDG^T)$. Both are centered, and for each
fixed $k\ge2$,
\begin{equation}
 |\Ee Y^k-\Ee V^k|\le K_{k,M}/N.
 \label{eq:weingarten-moment-comparison}
\end{equation}
\end{lemma}
\begin{proof}
The idea is to expand both moments and compare their coefficients.
First, for a real matrix $E$, define
\begin{equation}
 \mathcal T_p(E)=\sum_{i_1,\ldots,i_{2k}}
       \delta_p(i)\prod_{h=1}^k E_{i_{2h-1},i_{2h}}.
 \label{eq:weingarten-contraction}
\end{equation}
To interpret this sum, superimpose the pairing $p$ with $\epsilon=\{\{1,2\},\ldots,\{2k-1,2k\}\}$, whose edges
correspond to the matrix entries. The resulting components alternate
between pairing edges and matrix edges. Following a component with
$\ell$ matrix edges gives a trace
$\Tr(E^{\eta_1}\cdots E^{\eta_\ell})$, where
$\eta_j\in\{1,T\}$ records the direction in which the corresponding
entry is traversed. The index sums for different components are
separate, so $\mathcal T_p(E)$ is the product of these traces.

Now suppose $E$ is traceless. A component with only one matrix edge
contributes $\Tr E=0$, so any contraction containing such a component
vanishes. Every component of a nonzero contraction must therefore use
at least two matrix edges. There are $k$ matrix edges in total, which
means that at most $\lfloor k/2\rfloor$ components can remain.

For a word of length $\ell\ge2$, Frobenius Cauchy--Schwarz gives
\[
 |\Tr(C^{\eta_1}\cdots C^{\eta_\ell})|
 \le\|C\|_\F^2\|C\|_\op^{\ell-2}\le N^{\ell/2},\qquad
 |\Tr(D^{\eta_1}\cdots D^{\eta_\ell})|\le NM^\ell.
\]
The word lengths add up to $k$. Multiplying these bounds over the
components, and using the bound on their number, gives
\begin{equation}
 |\mathcal T_p(C^T)|\le N^{k/2},\qquad
 |\mathcal T_q(D)|\le M^kN^{\lfloor k/2\rfloor}.
 \label{eq:contraction-bounds}
\end{equation}
We can now return to the moments. Expand
$Y=\sum_{a,b,u,v}C_{ab}D_{uv}O_{bu}O_{av}$ and apply
Lemma~\ref{lem:weingarten-fixed-order}. For $V$, use Isserlis' formula
\cite{Isserlis} instead. This gives
\begin{align*}
 \Ee Y^k&=\sum_{p,q}\operatorname{Wg}_{N,k}(p,q)
                      \mathcal T_p(C^T)\mathcal T_q(D),\\
 \Ee V^k&=N^{-k}\sum_p\mathcal T_p(C^T)\mathcal T_p(D).
\end{align*}
Notice the difference between the two lines. A Gaussian covariance
identifies both the row indices and the column indices, so the same
pairing appears in both contractions in the second line. Subtracting
the two expansions, and applying the Weingarten estimate together with
\eqref{eq:contraction-bounds}, yields
\[
 |\Ee Y^k-\Ee V^k|
 \le L_k^2w_kM^kN^{-k-1+k/2+\lfloor k/2\rfloor}
 \le L_k^2w_kM^k/N.
\]
Lastly, the second-entry moment formulas for Haar and Gaussian matrices,
together with $\Tr C=\Tr D=0$, show that both means are zero. The bound
above depends only on $k$ and $M$, not on the particular matrices, as
required.
\end{proof}

We have now reduced the remaining moment estimates to a Gaussian
quadratic form. For a symmetric matrix $H$ and a standard Gaussian
vector $x$, the usual cumulant identity is
\begin{equation}
 \kappa_r(x^THx-\Tr H)=2^{r-1}(r-1)!\Tr(H^r),\qquad r\ge2.
 \label{eq:quadratic-cumulant}
\end{equation}
See \cite{Lancaster}. To obtain it, diagonalize $H$ and expand the
logarithm of the moment generating function,
$\sum_j[-t\lambda_j-\tfrac12\log(1-2t\lambda_j)]$.
The coefficient of each power of $t$ gives the corresponding cumulant.

\begin{theorem}\label{thm:weingarten-cumulants}
Under~\eqref{eq:weingarten-CD}, for every fixed $r\ge3$,
\begin{equation}
 |\kappa_r(\Tr(CODO^T))|
 \le K_{r,M}\bigl(N^{-1}+N^{1-r/2}\bigr).
 \label{eq:weingarten-cumulant-bound}
\end{equation}
\end{theorem}
\begin{proof}
We first work with the Gaussian statistic $V$. Stack the columns of
$G$ into one vector and let $x=\sqrt N\operatorname{vec}(G)$. Its
coordinates are independent standard Gaussian variables. The statistic
can then be written as $V=x^THx$, where
\[
 H=\frac{D^T\otimes C+D\otimes C^T}{2N},\qquad \Tr H=0.
\]
The two terms in the numerator make $H$ symmetric. Using the
Kronecker-product norm identities and \eqref{eq:weingarten-CD}, we get
\begin{equation}
 \|H\|_\op\le M/\sqrt N,\qquad \|H\|_\F\le M.
 \label{eq:weingarten-H-bounds}
\end{equation}
Thus the quadratic-form identity \eqref{eq:quadratic-cumulant} gives,
for $j\ge2$,
\begin{equation}
 |\kappa_j(V)|
 \le2^{j-1}(j-1)!\|H\|_\op^{j-2}\|H\|_\F^2
 \le2^{j-1}(j-1)!M^jN^{1-j/2}.
 \label{eq:weingarten-gaussian-cumulant-bound}
\end{equation}
Here the first inequality bounds two factors of $H$ in Frobenius norm
and the remaining $j-2$ factors in operator norm.

Now fix the order $r$. The cumulants just estimated are uniformly
bounded through order $r$. By \eqref{eq:moment-cumulant}, the moments
through that order are uniformly bounded as well.
Lemma~\ref{lem:weingarten-moment-comparison} says that each moment of
$Y$ differs from the corresponding moment of $V$ by
$O_{r,M}(N^{-1})$. In particular, the moments of $Y$ satisfy the same
boundedness property.

Recall that the $r$th cumulant is a fixed polynomial in the first $r$
moments. On a bounded set this polynomial is Lipschitz. Consequently,
$|\kappa_r(Y)-\kappa_r(V)|=O_{r,M}(N^{-1})$. Combining this error
with \eqref{eq:weingarten-gaussian-cumulant-bound} proves the claimed
bound. We have used only the assumptions in \eqref{eq:weingarten-CD}.
\end{proof}

\begin{corollary}
\label{cor:weingarten-uniform-cumulant}
Let $\widetilde A=A-\operatorname{tr}_N(A)I_N$ and
$\widetilde B=B-\operatorname{tr}_N(B)I_N$. For each fixed $r\ge3$,
\begin{equation}
 \sup_{B\in\mathrm O(N)}\ \sup_{\|A\|_\F^2=N}
 |\kappa_r(\Tr(\widetilde A O\widetilde B O^T))|
 \le K_r(N^{-1}+N^{1-r/2})\longrightarrow0.
 \label{eq:weingarten-uniform-cumulant}
\end{equation}
No trace conditions on $B$ are needed for this estimate. There is also a
constant $K$, independent of $N,A,B$, such that for all sufficiently large $N$,
\begin{equation}
 \sup_{B\in\mathrm O(N)}\ \sup_{\|A\|_\F^2=N}
 \Ee_O\!\left[\Tr(\widetilde A O\widetilde B O^T)^4\right]\le K.
 \label{eq:uniform-fourth-moment}
\end{equation}
\end{corollary}
\begin{proof}
First notice that $\widetilde A$ and $\widetilde B$ are traceless.
Also, $\|\widetilde A\|_\F^2=N(1-(\operatorname{tr}_N A)^2)\le N$
and $\|\widetilde B\|_\op\le2$. Thus the preceding theorem applies
with $M=2$, which proves the cumulant bound.

For the fourth moment, let
$Y=\Tr(\widetilde A O\widetilde B O^T)$. Subtracting the scalar parts
centers the statistic without changing its variance, so $Y$ is centered
and has the same variance as $\Tr(AOBO^T)$. Recall the quantities $x,y$
from the proof of Proposition~\ref{prop:variance-error}. We have
$x+y\le1$, $|\tau|\le1$, and $|b|\le1$. The coefficients of $x$ and
$y$ in the variance are at most $2$ and $4$, respectively. Hence
$\Ee_OY^2\le4$ for $N\ge2$.

Since $Y$ is centered, the moment--cumulant identity reads
$\Ee_OY^4=\kappa_4(Y)+3(\Ee_OY^2)^2\le|\kappa_4(Y)|+48$.
The fourth cumulant is uniformly bounded by the first part of the
corollary. This proves \eqref{eq:uniform-fourth-moment}.
\end{proof}

\begin{theorem}\label{thm:deterministic-clt}
Let $B_N\in\mathrm O(N)$ satisfy
$\operatorname{tr}_N B_N\to0$ and $\operatorname{tr}_N(B_N^2)\to0$.
If $\|A_N\|_\F^2=N$ and $a_N=\operatorname{tr}_N A_N\to s$, then
\[
 \Tr(A_NO_NB_NO_N^T)-a_N\Tr B_N\dto N(0,1-s^2).
\]
This holds along every deterministic dimension sequence tending to infinity.
\end{theorem}
\begin{proof}
The centered statistic is $\Tr(\widetilde A_NO_N\widetilde B_NO_N^T)$,
after subtracting the scalar parts of $A_N$ and $B_N$. Its mean is zero. Corollary~\ref{cor:variance-limit} shows that its variance tends
to $1-s^2$, and Corollary~\ref{cor:weingarten-uniform-cumulant} shows
that each fixed cumulant of order at least three tends to zero. These
are exactly the conditions in the moment criterion described above.
Thus the limit is $N(0,1-s^2)$. When $|s|=1$, the variance tends to zero,
and the second moment alone gives convergence to zero. Finally, all the
estimates depend on the dimension rather than on its index, so the same
proof applies along any deterministic dimension sequence tending to
infinity.
\end{proof}

The coefficient estimate in
Lemma~\ref{lem:weingarten-fixed-order} is a coarse fixed-order
consequence of Collins--Matsumoto
\cite[Theorem~4.11]{CollinsMatsumoto}. What lets us use only a
Frobenius bound on the test is the uniform estimate for the
contractions. For comparison, Chatterjee--Meckes
\cite[Theorems~4.4--4.5]{ChatterjeeMeckes} and Meckes
\cite[Chapter~2]{Meckes} give quantitative normal approximation for
linear projections $\Tr(AO)$ of Haar matrices. Here we prove
convergence for $\Tr(AOBO^T)$, but do not turn the fixed-order
cumulant estimates into a bound on the distributional error. Obtaining
a rate would require an additional approximation argument, as well
as rates for the trace and test-normalization limits.

Now suppose the conjugated matrix is random. Applying the preceding
theorem separately to deterministic sequences is not by itself the
conditional statement we need. We must control its characteristic
function uniformly on a set of matrices whose probability tends to
one. Once we have this control, we can also prove independence from
any convergent observable of the underlying randomness.

\begin{theorem}[Conditional Haar-conjugation central limit theorem]\label{thm:linear}
Let $Q_n\in\mathrm O(d_n)$ be measurable with respect to a sigma-field
$\mathcal F_n$, where $d_n\to\infty$, and let $O_n$ be Haar and independent
of $\mathcal F_n$. Suppose only that
\begin{equation}
 \frac{\Tr Q_n}{d_n}\longrightarrow0,\qquad
 \frac{\Tr(Q_n^2)}{d_n}\longrightarrow0
 \quad\text{in probability}.
 \tag{C}\label{eq:character-condition}
\end{equation}
For deterministic real tests satisfying~\eqref{eq:test-normalization}, put
$X_n=\Tr(A_nO_nQ_nO_n^T)-s_n\Tr Q_n$. Then $|s|\le1$, and for every
$t\in\R$,
\begin{equation}
 \Ee[e^{itX_n}\mid\mathcal F_n]
 \longrightarrow c_t:=\exp\!\left(-\frac{t^2(1-s^2)}2\right)
 \quad\text{in probability and in }L^1.
 \label{eq:conditional-transfer}
\end{equation}
In fact, for $\mathcal F_n$-measurable $Z_n$,
\begin{equation}
 \sup_{ |Z_n|\le1}
 \left|\Ee[Z_ne^{itX_n}]-c_t\Ee Z_n\right|\longrightarrow0.
 \label{eq:bounded-environment-transfer}
\end{equation}
If a fixed-dimensional $\mathcal F_n$-measurable real vector $W_n$ satisfies
$W_n\dto W$, then
$(W_n,X_n)\dto(W,\sqrt{1-s^2}G)$, where $G\sim N(0,1)$ is independent
of $W$. No tightness or convergence of $\Tr Q_n$ is assumed.
\end{theorem}

The two conditions in \eqref{eq:character-condition} are necessary
and sufficient for the conditional variance limit when that limit
is required uniformly over normalized tests. For one particular
test, cancellation may make a condition unnecessary.
Remark~\ref{rem:trace-sharpness} states the distinction and gives
an example.

\begin{proof}
We first assume exact normalization, $\|A_n\|_\F^2=d_n$. The goal is
to apply the deterministic theorem to the possible values of $Q_n$.
By \textup{(C)}, we can choose $\delta_n\downarrow0$ such that
$\Pp((|\Tr Q_n|+|\Tr Q_n^2|)/d_n>\delta_n)\to0$. Define
\[
 \mathcal G_n=\{B\in\mathrm O(d_n):
          (|\Tr B|+|\Tr B^2|)/d_n\le\delta_n\}.
\]
These sets are compact, and $\Pp(Q_n\in\mathcal G_n)\to1$. In
particular, they are nonempty for all sufficiently large $n$.
For a deterministic matrix $B$, write
\[
 \phi_n(B;t)=\Ee_{O_n}\exp\!\left(
 it[\Tr(A_nO_nBO_n^T)-s_n\Tr B]\right).
\]
The integrand is bounded and jointly continuous in $(O_n,B)$. By
dominated convergence, $\phi_n(\cdot;t)$ is continuous, and hence
measurable. If we select any deterministic $B_n\in\mathcal G_n$, its
two normalized traces tend to zero by the definition of $\mathcal G_n$.
Theorem~\ref{thm:deterministic-clt} therefore gives
$\phi_n(B_n;t)\to c_t$. We claim that this convergence is uniform on
$\mathcal G_n$:
\begin{equation}
 \sup_{B\in\mathcal G_n}|\phi_n(B;t)-c_t|\longrightarrow0.
 \label{eq:uniform-cf-linear}
\end{equation}
To see this, suppose the claim were false. We could then find a
subsequence and a matrix from each corresponding $\mathcal G_n$ for
which the difference stays bounded away from zero. Choose arbitrary
matrices from the other nonempty sets to complete a deterministic
sequence. Its two normalized traces still tend to zero, contradicting
the deterministic theorem. This argument also allows nonconsecutive
or repeated dimensions, provided $d_n\to\infty$.

We next return to the random matrix. Since $O_n$ is independent of
$\mathcal F_n$, independence and Fubini give
$\Ee[e^{itX_n}\mid\mathcal F_n]=\phi_n(Q_n;t)$. Split the expectation
according to whether $Q_n$ belongs to $\mathcal G_n$. On that set we
use the supremum in \eqref{eq:uniform-cf-linear}; on its complement
we use the bound $2$ for the difference of characteristic functions.
Thus
\begin{equation}
 \Ee|\phi_n(Q_n;t)-c_t|
 \le\sup_{B\in\mathcal G_n}|\phi_n(B;t)-c_t|
      +2\Pp(Q_n\notin\mathcal G_n)\longrightarrow0.
 \label{eq:conditional-cf-bound-linear}
\end{equation}
This proves convergence in $L^1$, and therefore in probability, as
asserted in \eqref{eq:conditional-transfer}. If $Z_n$ is
$\mathcal F_n$-measurable and $|Z_n|\le1$, multiplication by $Z_n$
and the conditional-expectation identity bound the error in
\eqref{eq:bounded-environment-transfer} by the same $L^1$ quantity.
Taking $Z_n=e^{iu^TW_n}$ gives convergence of the joint characteristic
function to the product of the two limiting characteristic functions.
L\'evy's continuity theorem now proves both joint convergence and
independence.

It remains to allow approximate normalization. Let
$r_n=(d_n^{-1}\|A_n\|_\F^2)^{1/2}\to1$, and denote by
$\widehat X_n$ the centered statistic with test $A_n/r_n$. The case
already proved makes $\widehat X_n$ tight, and
$X_n=r_n\widehat X_n$. For fixed $t$, we have
$\Ee|e^{itr_n\widehat X_n}-e^{it\widehat X_n}|\to0$. Indeed, on
$\{|\widehat X_n|\le M\}$ the difference is at most
$|t|M|r_n-1|$, while tightness makes the contribution from the complement
small when $M$ is large. Conditional expectation does not increase the
$L^1$ norm, so this error also tends to zero after conditioning. All
conclusions therefore remain valid. Finally,
$|s_n|\le\|A_n\|_\F/\sqrt{d_n}\to1$ implies $|s|\le1$.
\end{proof}

\begin{corollary}\label{cor:companion-law}
\emergencystretch=1em
Under the hypotheses of Theorem~\ref{thm:linear}, suppose additionally that
$(W_n,\Tr Q_n)\dto(W,\Xi)$ for an $\mathcal F_n$-measurable vector $W_n$.
Then
\[
 (W_n,\Tr Q_n,X_n)\dto(W,\Xi,\sqrt{1-s^2}G),\qquad
 \Tr(A_nO_nQ_nO_n^T)\dto s\Xi+\sqrt{1-s^2}G,
\]
where $G$ is independent of $(W,\Xi)$. Convergence of $\Tr Q_n$ is
assumed here but not in Theorem~\ref{thm:linear}. In the hypotheses of
this corollary, it already implies the first normalized-trace condition
when $d_n\to\infty$.
\end{corollary}

\begin{proof}
Apply Theorem~\ref{thm:linear} to the vector $(W_n,\Tr Q_n)$.
This gives the joint limit with a Gaussian variable independent of
$(W,\Xi)$. To recover the uncentered statistic, we must add back
$s_n\Tr Q_n$. Since $\Tr Q_n$ converges in distribution, it is tight;
together with $s_n\to s$, this gives
$(s_n-s)\Tr Q_n\to0$ in probability. Slutsky's theorem then gives the
stated uncentered limit.
\end{proof}

We now return to $Q_n=\rho_n(\sigma_n)$. The representation structure
will let us relate its two normalized traces. Apart from the trivial
and sign representations, every irreducible representation of $S_n$
has dimension tending to infinity uniformly with $n$. We will use
this lower bound in a tensor-square calculation to control the
character at $\sigma_n^2$, without any further restriction on the
shapes of the constituent diagrams.

Recall that a partition $\lambda\vdash n$ is a decreasing sequence
of positive integers with sum $n$. Its Young diagram has
$\lambda_i$ left-justified boxes in row $i$. We use English notation,
so the first row is at the top. Write $\lambda'$ for the transposed
diagram, $S^\lambda$ for the real Specht module, and
$f^\lambda=\dim S^\lambda$. The real Specht modules are absolutely
irreducible and give all irreducible real representations of $S_n$
\cite{Sagan}. Averaging a positive inner product makes any realization
orthogonal.

The randomized law does not depend on which orthogonal realization
we choose. To see this, suppose $T$ intertwines two such realizations.
Then $T^TT$ commutes with the first one, and the polar factor
$U=T(T^TT)^{-1/2}$ is an orthogonal intertwiner. Right multiplication
by $U$ preserves the Haar law. Since it also preserves independence
from $\sigma_n$, the two randomized matrices have the same law.

Let us write the representation as a sum of its irreducible
constituents, keeping track of the multiplicities:
\begin{equation}
 \rho_n\cong\bigoplus_{\lambda\vdash n}(S^\lambda)^{\oplus m_{\lambda,n}},
 \qquad d_n=\sum_{\lambda\vdash n}m_{\lambda,n}f^\lambda.
 \label{eq:general-decomposition}
\end{equation}
Define
\begin{equation}
 t_n=m_{(n),n}+m_{(1^n),n},\qquad
 r_n=\sum_{\lambda\vdash n}m_{\lambda,n}.
 \label{eq:one-dimensional-mass}
\end{equation}
Thus $t_n$ counts the trivial and sign summands, while $r_n$ counts
all irreducible summands, in both cases with multiplicity. Since the
two exceptional representations are one-dimensional, $t_n/d_n$ is
their fraction of the total dimension. It is not their fraction of
the number of summands. The other multiplicities and diagrams may
vary arbitrarily with $n$.

For real class functions, we use the normalized inner product
$\langle f,g\rangle_{S_n}=\Ee[f(\sigma_n)g(\sigma_n)]$.
Irreducible characters are orthonormal for this inner product.
In particular,
\begin{equation}
 \Ee[\chi^{\lambda_n}(\sigma_n)^2]=1,\qquad
 \Pp(|\chi^{\lambda_n}(\sigma_n)|>M)\le M^{-2}.
 \label{eq:irreducible-tightness}
\end{equation}
The second inequality follows from the first by Markov's inequality.
It shows that every irreducible-character sequence is tight, but
tightness alone does not give a unique limiting law. We will also
use the following two classical facts to control the normalized
character at the square.

\begin{lemma}\label{lem:classical-character-inputs}
For every $\lambda\vdash n$,
\begin{equation}
 \Ee[\chi^\lambda(\sigma_n^2)]=1.
 \label{eq:frobenius-schur}
\end{equation}
For $n\ge5$ and $\lambda\notin\{(n),(1^n)\}$,
\begin{equation}
 f^\lambda\ge n-1.
 \label{eq:minimum-degree}
\end{equation}
\end{lemma}
\begin{proof}
Choose an orthogonal realization $R$ of $S^\lambda$ and write
$f=f^\lambda$ for its dimension. Since this real representation is
absolutely irreducible, matrix-coefficient orthogonality gives
$\Ee[R(\sigma_n)_{ij}R(\sigma_n)_{k\ell}]
=f^{-1}\delta_{ik}\delta_{j\ell}$. Now expand the trace of
$R(\sigma_n)^2$ entry by entry:
\[
 \Ee[\chi^\lambda(\sigma_n^2)]
 =\sum_{i,j=1}^f\Ee[R(\sigma_n)_{ij}R(\sigma_n)_{ji}]
 =\frac1f\sum_{i,j=1}^f\delta_{ij}=1.
\]
Only the terms with $i=j$ remain, which explains the value $1$.
This is the Frobenius--Schur identity in the real setting used here.

The second assertion is the classical minimum-degree bound of
Rasala \cite{Rasala}; see also
Bessenrodt--Nguyen--Olsson--Tong-Viet \cite[Lemma~3.1(1)]{BNOT}
for the range $n\ge5$. We only need this lower bound, not a
classification of the partitions for which equality holds.
\end{proof}

\begin{proposition}
\label{prop:square-character-bound}
Let $\chi$ be the character of any real representation $V$ of $S_n$, of
dimension $d>0$, and let $\sigma$ be uniform on $S_n$. For $n\ge5$,
\begin{equation}
 \Ee\!\left[\left(\frac{\chi(\sigma^2)}d\right)^2\right]
 \le 2\Ee\!\left[\left(\frac{\chi(\sigma)}d\right)^2\right]
       +\frac1{n-1}.
 \label{eq:square-character-bound}
\end{equation}
\end{proposition}
\begin{proof}
The useful observation is that $\chi^2$ is itself the character of a
representation, namely $V\otimes V$. We can therefore decompose this
representation with nonnegative integer multiplicities:
\[
 V\otimes V\cong\bigoplus_{\nu\vdash n}(S^\nu)^{\oplus c_\nu},
 \qquad c_\nu\in\mathbb Z_{\ge0},\qquad
 \sum_{\nu\vdash n}c_\nu f^\nu=d^2.
\]
Now evaluate the character at $\sigma^2$ and take expectations.
Each irreducible constituent has expected character $1$ at the square
by \eqref{eq:frobenius-schur}. Thus
\begin{equation}
 \Ee[\chi(\sigma^2)^2]=\sum_{\nu\vdash n}c_\nu.
 \label{eq:tensor-square-multiplicity}
\end{equation}
We have reduced the problem to bounding the number of constituents,
counted with multiplicity. First separate the two one-dimensional
representations. Character orthogonality gives
\[
 c_{(n)}=\langle\chi^2,1\rangle_{S_n}=\Ee[\chi(\sigma)^2],\qquad
 0\le c_{(1^n)}=\Ee[\sgn(\sigma)\chi(\sigma)^2]
                   \le\Ee[\chi(\sigma)^2].
\]
For the sign constituent, nonnegativity comes from its interpretation
as a multiplicity. The upper bound follows because
$\sgn(\sigma)\le1$ and $\chi(\sigma)^2\ge0$.

Every other constituent has dimension at least $n-1$, by
\eqref{eq:minimum-degree}. Consequently
\[
 \sum_{\nu\notin\{(n),(1^n)\}}c_\nu
 \le\frac1{n-1}\sum_{\nu\notin\{(n),(1^n)\}}c_\nu f^\nu
 \le\frac{d^2}{n-1}.
\]
The last inequality uses the total dimension $d^2$ of $V\otimes V$.
Add the bounds for the two one-dimensional constituents and the
remaining constituents in \eqref{eq:tensor-square-multiplicity}, then
divide by $d^2$. This proves the result without having to compute the
individual multiplicities $c_\nu$.
\end{proof}

It is worth distinguishing the two uses of a square here. The
function $g\mapsto\chi(g^2)$ is the virtual character
$\psi^2\chi=\chi_{\operatorname{Sym}^2V}-\chi_{\Lambda^2V}$:
subtracting the two eigenvalue formulas leaves $\sum_i z_i^2$.
Consequently,
\begin{equation}
 \Ee[\chi(\sigma_n^2)^2]
 =\langle\psi^2\chi,\psi^2\chi\rangle_{S_n}.
 \label{eq:square-character-norm}
\end{equation}
This norm need not be one, even when $V$ is irreducible. In
Proposition~\ref{prop:square-character-bound}, we instead decomposed
the character $\chi^2$ of the genuine representation $V\otimes V$
and then evaluated it at $\sigma_n^2$. This order of operations
lets us use nonnegative multiplicities, rather than decompose the
virtual character $\psi^2\chi$ itself.

\begin{proposition}
\label{prop:multiplicity-traces}
For the decomposition~\eqref{eq:general-decomposition}, put
$\beta_n=\chi_n(\sigma_n)/d_n$ and
$\tau_n=\chi_n(\sigma_n^2)/d_n$. For $n\ge5$,
\begin{align}
 \frac12\left(\frac{t_n}{d_n}\right)^2
 &\le\Ee[\beta_n^2]
 =\frac{\sum_{\lambda\vdash n}m_{\lambda,n}^2}{d_n^2}
 \le\left(\frac{t_n}{d_n}\right)^2+\frac1{(n-1)^2},
 \label{eq:ordinary-multiplicity-bound}\\
 \Ee[\tau_n]&=\frac{r_n}{d_n},\qquad
 \frac{t_n}{d_n}\le\frac{r_n}{d_n}
       \le\frac{t_n}{d_n}+\frac1{n-1}.
 \label{eq:square-mean-multiplicity}
\end{align}
As $n\to\infty$, the following are equivalent:
\begin{equation}
 \frac{t_n}{d_n}\longrightarrow0
 \quad\Longleftrightarrow\quad
 \beta_n\longrightarrow0\ \text{in probability}
 \quad\Longleftrightarrow\quad
 \tau_n\longrightarrow0\ \text{in probability}.
 \label{eq:trace-mass-equivalence}
\end{equation}
Each probabilistic convergence may equivalently be required in $L^2$.
The condition $r_n/d_n\to0$ is also equivalent to these conditions.
\end{proposition}
\begin{proof}
First apply character orthogonality to the decomposition of $\chi_n$.
All cross terms vanish, giving the equality in
\eqref{eq:ordinary-multiplicity-bound}. The two one-dimensional
multiplicities have sum $t_n$, so the sum of their squares lies between
$t_n^2/2$ and $t_n^2$. For every other constituent, the dimension is at
least $n-1$. Since these constituents account for $d_n-t_n$ dimensions,
we get
\[
 \sum_{\lambda\notin\{(n),(1^n)\}}m_{\lambda,n}
 \le\frac{d_n-t_n}{n-1}.
\]
Their multiplicities are nonnegative, so the sum of their squares is
at most the square of their sum. Combining this observation with the
bounds for the two one-dimensional multiplicities proves
\eqref{eq:ordinary-multiplicity-bound}.

The left side of the last display is $r_n-t_n$, which also gives the
two bounds on $r_n/d_n$. Further, each irreducible character has
expectation $1$ at $\sigma_n^2$. Summing \eqref{eq:frobenius-schur}
with the multiplicities therefore gives $\Ee[\tau_n]=r_n/d_n$.
This proves \eqref{eq:square-mean-multiplicity}.

We now prove the convergence assertions. Since $\rho_n(\sigma)$ is
orthogonal, $|\chi_n(\sigma)|\le d_n$. Thus
$|\beta_n|,|\tau_n|\le1$. For either sequence, this bound makes
convergence in probability to zero equivalent to convergence in $L^2$.
For example, for every $\delta>0$,
$\Ee[\beta_n^2]\le\delta^2+\Pp(|\beta_n|>\delta)$.
First let $n\to\infty$, then let $\delta\downarrow0$. The reverse
implication follows from Markov's inequality.

The two bounds in \eqref{eq:ordinary-multiplicity-bound} now show that
$t_n/d_n\to0$ if and only if $\beta_n\to0$.
Proposition~\ref{prop:square-character-bound} then gives
$\tau_n\to0$. Conversely, suppose $\tau_n\to0$ in probability.
Boundedness gives $\Ee[\tau_n]\to0$, and
$\Ee[\tau_n]=r_n/d_n\ge t_n/d_n\ge0$. Hence $t_n/d_n\to0$ as well.
Finally, the bounds in \eqref{eq:square-mean-multiplicity} show that
$r_n/d_n\to0$ is equivalent to $t_n/d_n\to0$.
\end{proof}

We have now replaced the two normalized-trace conditions by one
condition on multiplicities. The next step is to show that this
condition is necessary for the distributional conclusion for all
tests. Comparing candidate variances is not enough: convergence in
distribution need not give convergence of second moments. The
uniform fourth-moment bound will provide the uniform integrability
needed for that step.

\begin{theorem}[Universal centered limits]\label{thm:universality}
Let $n\to\infty$, let $\rho_n:S_n\to\mathrm O(d_n)$ be arbitrary real
representations with $d_n\to\infty$, and use
\eqref{eq:general-decomposition}--\eqref{eq:one-dimensional-mass}.
Let $\sigma_n$ be uniform on $S_n$ and $O_n$ be independent Haar orthogonal.
For a deterministic test satisfying~\eqref{eq:test-normalization}, put
\[
 X_n=\Tr(A_nO_n\rho_n(\sigma_n)O_n^T)-s_n\chi_n(\sigma_n).
\]
Then $t_n/d_n\to0$ if and only if, for \emph{every} such test sequence
and every $u\in\R$,
\begin{equation}
 \Ee[e^{iuX_n}\mid\sigma_n]
 \longrightarrow\exp\!\left(-\frac{u^2(1-s^2)}2\right)
 \quad\text{in probability and in }L^1.
 \label{eq:universal-conditional-limit}
\end{equation}
The same condition is equivalent to requiring only
$X_n\dto N(0,1-s^2)$ for every such deterministic test sequence.
Under this condition the bounded-environment estimate
\eqref{eq:bounded-environment-transfer} holds with
$\mathcal F_n=\sigma(\sigma_n)$, and any convergent finite-dimensional
$\sigma_n$-measurable vector is jointly asymptotically independent of the
Gaussian fluctuation. No convergence or tightness of the ordinary
character is required.
\end{theorem}
\begin{proof}
Suppose first that $t_n/d_n\to0$. By
Proposition~\ref{prop:multiplicity-traces}, both $\beta_n$ and
$\tau_n$ tend to zero in $L^2$. Recall that
$Q_n^2=\rho_n(\sigma_n^2)$, so these are precisely the two normalized
traces in \textup{(C)}. We may therefore apply Theorem~\ref{thm:linear}
with $\mathcal F_n=\sigma(\sigma_n)$. This gives
\eqref{eq:universal-conditional-limit}, the bounded-environment estimate,
and the joint limit with every convergent vector specified in the
theorem. Taking expectations of the conditional characteristic
functions also gives the unconditional Gaussian limit.

For the converse, it is enough to use one test. Take
\[
 A_n^-=\sqrt{d_n/2}\,(E_{12}-E_{21}).
\]
Since $d_n\to\infty$, this matrix is defined for all sufficiently large
$n$. It is traceless and has squared Frobenius norm $d_n$. Write
$X_n^-=\Tr(A_n^-O_nQ_nO_n^T)$. Its conditional mean is zero, so its
conditional second moment is its conditional variance. The exact
variance formula gives
\begin{equation}
 \Ee[(X_n^-)^2\mid\sigma_n]
 =\frac{d_n}{d_n-1}(1-\tau_n),\qquad
 \Ee[(X_n^-)^2]=\frac{d_n-r_n}{d_n-1}.
 \label{eq:skew-diagnostic-multiplicity}
\end{equation}
The second equality follows by taking expectations and using
\eqref{eq:square-mean-multiplicity}.

We want to deduce convergence of these second moments from the assumed
convergence in distribution. For this, we need uniform integrability.
Apply \eqref{eq:uniform-fourth-moment} conditionally to each value of
$Q_n$, then take expectations. After discarding finitely many terms,
$\sup_n\Ee[(X_n^-)^4]<\infty$. In particular, the expectation of
$(X_n^-)^2$ on $\{|X_n^-|>M\}$ is at most
$M^{-2}\sup_n\Ee[(X_n^-)^4]$, which tends to zero as $M\to\infty$.
Thus the squares are uniformly integrable.

The unconditional conclusion for this test is $X_n^-\dto N(0,1)$.
Uniform integrability now gives $\Ee[(X_n^-)^2]\to1$.
Since $d_n\to\infty$, \eqref{eq:skew-diagnostic-multiplicity} implies
$r_n/d_n\to0$, and therefore $t_n/d_n\to0$. We have proved necessity
from the unconditional conclusion alone. This completes both
equivalences.
\end{proof}

Notice the quantifier over tests in the theorem. We ask for the
specified variance $1-s^2$ for every normalized test, not just for
some Gaussian limit of one specially chosen statistic. In the
symmetric-group setting, the skew-symmetric test used in the proof
already detects failure of the condition. This is consistent with
Remark~\ref{rem:trace-sharpness}. For arbitrary orthogonal matrices,
the two normalized traces need not be related. Here
Proposition~\ref{prop:multiplicity-traces} uses both the uniformity
of $\sigma_n$ and the fact that the constituents are representations
of $S_n$.

We can now add back the character without making a separate
assumption at the square. If the ordinary character is tight,
dividing it by $d_n\to\infty$ makes the first normalized character
small. Proposition~\ref{prop:square-character-bound} then makes the
second one small as well. This gives a necessary-and-sufficient
condition for the uncentered limits.

\begin{theorem}[Character convergence and uncentered limits]
\label{thm:character-equivalence}
Let $n\to\infty$, let $\rho_n:S_n\to\mathrm O(d_n)$ be arbitrary real
representations with $d_n\to\infty$, let $\sigma_n$ be uniform, and let
$O_n$ be independent Haar orthogonal. No multiplicity hypothesis is imposed.
For a real random variable $\Xi$ with a fixed law, the convergence
$\chi_n(\sigma_n)\dto\Xi$ holds if and only if, for every deterministic test
satisfying~\eqref{eq:test-normalization},
\begin{equation}
 \Tr(A_nO_n\rho_n(\sigma_n)O_n^T)
 \dto s\Xi+\sqrt{1-s^2}\,G,
 \qquad G\sim N(0,1),\quad G\indep\Xi.
 \label{eq:character-equivalence-limit}
\end{equation}
When these assertions hold, $t_n/d_n\to0$ and
$\chi_n(\sigma_n^2)/d_n\to0$ in $L^2$ automatically. More precisely, if
$(W_n,\chi_n(\sigma_n))\dto(W,\Xi)$ for a fixed-dimensional
$\sigma_n$-measurable vector $W_n$, then
\begin{equation}
 (W_n,\chi_n(\sigma_n),X_n)
 \dto(W,\Xi,\sqrt{1-s^2}\,G),
 \qquad G\indep(W,\Xi).
 \label{eq:character-joint-upgrade}
\end{equation}
\end{theorem}
\begin{proof}
Suppose first that $\chi_n(\sigma_n)\dto\Xi$. The character sequence
is then tight. Since $d_n\to\infty$, dividing by $d_n$ gives
$\beta_n=\chi_n(\sigma_n)/d_n\to0$ in probability. We also have
$|\beta_n|\le1$, so the convergence holds in $L^2$.
Proposition~\ref{prop:square-character-bound} gives $\tau_n\to0$ in
$L^2$, and Proposition~\ref{prop:multiplicity-traces} gives
$t_n/d_n\to0$.

Thus Theorem~\ref{thm:universality} applies. Use it with the convergent
vector $(W_n,\chi_n(\sigma_n))$ to obtain
\eqref{eq:character-joint-upgrade}. To recover the uncentered statistic,
add back $s_n\chi_n(\sigma_n)$. Tightness and $s_n\to s$ imply
$(s_n-s)\chi_n(\sigma_n)\to0$ in probability, so Slutsky's theorem
gives \eqref{eq:character-equivalence-limit}. Notice that the argument
never assumes convergence of the unnormalized second moments of the
characters.

For the reverse implication, take $A_n=I_{d_n}$. Then $s_n=1$, the
Gaussian term vanishes, and the statistic is exactly
$\chi_n(\sigma_n)$. The conclusion for this test is therefore
$\chi_n(\sigma_n)\dto\Xi$, as desired.
\end{proof}

The main use of the representation structure is in the forward
implication: ordinary character convergence gives both normalized
trace conditions. The reverse implication just uses the identity
test, which returns the character itself. Thus the remaining
uncentered problem is to determine whether the ordinary character
converges and, when it does, identify its law. This equivalence does
not classify all Young-diagram sequences with a character limit.

\begin{corollary}\label{thm:entries}
Under $t_n/d_n\to0$, put $M_n=O_n\rho_n(\sigma_n)O_n^T$. Fix $r\ge1$
and, for each $n$, choose pairwise distinct ordered pairs
$(a_j^{(n)},b_j^{(n)})\in\{1,\ldots,d_n\}^2$. Then
\begin{equation}
 \left(\sqrt{d_n}\left((M_n)_{a_j^{(n)}b_j^{(n)}}
 -\frac{\chi_n(\sigma_n)}{d_n}
     \mathbf1_{\{a_j^{(n)}=b_j^{(n)}\}}\right)\right)_{j=1}^r
 \dto N(0,I_r).
 \label{eq:centered-entry-universality}
\end{equation}
The limit is jointly independent of any convergent fixed-dimensional
$\sigma_n$-measurable vector. The centering can be omitted whenever
$\chi_n(\sigma_n)/\sqrt{d_n}\to0$ in probability; it is identically zero
for off-diagonal positions. In particular, ordinary character convergence
is sufficient to omit it.
\end{corollary}
\begin{proof}
We use Cram\'er--Wold, so fix $\theta\in\R^r\setminus\{0\}$ and
consider the test
\[
 A_n(\theta)=\frac{\sqrt{d_n}}{\|\theta\|_2}
                   \sum_{j=1}^r\theta_jE_{b_j^{(n)}a_j^{(n)}}.
\]
The positions are distinct, hence the matrix units in this sum are
Frobenius-orthonormal. It follows that $\|A_n(\theta)\|_\F^2=d_n$.
Only the diagonal positions contribute to the trace. By
Cauchy--Schwarz,
\[
 \left|\frac{\Tr A_n(\theta)}{d_n}\right|
 =\left|\frac{\sum_j\theta_j\mathbf1_{\{a_j^{(n)}=b_j^{(n)}\}}}
                  {\sqrt{d_n}\|\theta\|_2}\right|
 \le\sqrt{r/d_n}\longrightarrow0.
\]
Now recall that $\Tr(E_{ba}M_n)=(M_n)_{ab}$. Thus the centered
statistic for $A_n(\theta)$ is $\|\theta\|_2^{-1}$ times the linear
combination with coefficients $\theta_j$ of the entries in
\eqref{eq:centered-entry-universality}. Theorem~\ref{thm:universality}
gives a standard normal limit for this normalized combination.
Cram\'er--Wold therefore gives the vector limit $N(0,I_r)$.

For the joint assertion, apply the bounded-environment estimate with
$e^{iv^TW_n}$. The same characteristic-function argument as before
shows that the Gaussian vector is independent of the limiting $W$.
Finally, the extra hypothesis
$\chi_n(\sigma_n)/\sqrt{d_n}\to0$ makes each diagonal centering term
tend to zero in probability. There are only finitely many coordinates,
so Slutsky's theorem removes them together. Off the diagonal there is
no centering term to remove. Ordinary character convergence implies
tightness, which is enough for the extra hypothesis since
$d_n\to\infty$.
\end{proof}

\begin{corollary}
\label{cor:all-irreducible-entries}
Let $\lambda_n\vdash n$ be arbitrary, with
$d_n=f^{\lambda_n}\to\infty$, and let $\rho_n$ be any orthogonal
realization of $S^{\lambda_n}$. The centered conditional limit of
Theorem~\ref{thm:universality} holds, and every fixed collection of distinct
ordered entries satisfies
\begin{equation}
 \left(\sqrt{d_n}\,(O_n\rho_n(\sigma_n)O_n^T)_{a_j^{(n)}b_j^{(n)}}
       \right)_{j=1}^r\dto N(0,I_r).
 \label{eq:all-irreducible-entries}
\end{equation}
No stabilization or geometric restriction on $\lambda_n$ is required,
and the ordinary characters need not converge in distribution.
\end{corollary}
\begin{proof}
Since $d_n\to\infty$, the representation is eventually neither trivial
nor sign. Hence $t_n=0$, and Theorem~\ref{thm:universality} applies.
We can also see the normalized square-character convergence directly:
combine \eqref{eq:irreducible-tightness} with
Proposition~\ref{prop:square-character-bound} to get
\[
 \Ee\!\left[\left(\frac{\chi^{\lambda_n}(\sigma_n^2)}{d_n}\right)^2\right]
 \le\frac{2}{d_n^2}+\frac1{n-1}\longrightarrow0.
\]
For the centering, character orthogonality gives
$\Ee[(\chi^{\lambda_n}(\sigma_n)/\sqrt{d_n})^2]=1/d_n\to0$.
Thus the extra hypothesis of Corollary~\ref{thm:entries} holds, and
that corollary gives the uncentered entry limit.
\end{proof}

\begin{remark}\label{rem:centered-entries}
For a reducible representation, $t_n=o(d_n)$ need not imply
$\chi_n(\sigma_n)/\sqrt{d_n}\to0$. Thus the centered entry result
applies more generally than the uncentered one. There is a related
point about the scalar part of a test: $s_n\to0$ does not by itself
make $s_n\chi_n(\sigma_n)$ tend to zero if the character is not tight.
It is the product that must vanish. An exactly traceless test, on
the other hand, removes this term at every $n$. The examples in
Section~\ref{sec:characters} will distinguish these cases.
\end{remark}

\newpage

\section{Explicit character laws}\label{sec:characters}
\label{sec:diagram-polynomials}\label{sec:permutations}
\label{sec:examples}\label{sec:worked-tail}

We have determined when the centered Gaussian limit holds for arbitrary
real representations, and we have obtained entry limits for every
growing-dimensional irreducible sequence. We now turn to a different
question: what is the limiting law $\Xi$ of the ordinary character?
This law still contributes when the test has a scalar component. Once
we know it, Theorem~\ref{thm:character-equivalence} gives the
uncentered randomized-basis limit, without another assumption on the
normalized character at the square.

We consider fixed first-row tails, their transposes, and diagrams
with sublinear largest hook, together with suitable direct sums.
These families give explicit character laws, but they do not exhaust
all possible character limits. Exact eventual cycle-polynomial
identities are more restrictive: they satisfy the converse in
Proposition~\ref{prop:exact-tail-equivalence} and
Corollary~\ref{cor:row-column-polynomials}. We will also calculate
unnormalized characters at the square. These calculations give more
information than normalized convergence to zero, even though the
latter already follows from the general theorem.

Write $\ell(\lambda)=\lambda'_1$ for the number of rows, and let
\begin{equation}
 h(\lambda)=\lambda_1+\lambda'_1-1.
 \label{eq:largest-hook}
\end{equation}
This is the hook length of the upper-left box, which is the largest
hook length in the diagram. When this hook is small compared with
$n$, a long cycle will force the ordinary character to be exactly
zero. The same event will work for all constituents, even when their
multiplicities are unrestricted. We first treat fixed tails, where
we can instead write down the character explicitly.

Suppose $\mu\vdash k$ is fixed. We keep this diagram below the first
row and let the first row grow. More precisely, set
$\mu[n]=(n-k,\mu_1,\ldots,\mu_{\ell(\mu)})$ whenever
$n-k\ge\mu_1$, so that this is a partition. We also use the convention
$\varnothing[n]=(n)$. Thus $k$ counts the boxes below the first row,
while $\mu$ specifies their shape. Fixing their number alone would
not fix the tail.

We want a formula whose size depends on the fixed tail, rather than
on $n$. For $\alpha\vdash r$, let $m_j(\alpha)$ be the number of
parts of $\alpha$ equal to $j$. Define
\begin{align}
 F_\nu(x)&=\sum_{\alpha\,\vdash|\nu|}\chi^\nu(\alpha)
                 \prod_{j\ge1}\binom{x_j}{m_j(\alpha)},
                 \qquad F_\varnothing=1,\label{eq:induced-polynomial}\\
 q_\mu(x)&=\sum_{\substack{\nu\subseteq\mu\\
                    \mu/\nu\text{ a vertical strip}}}
           (-1)^{|\mu|-|\nu|}F_\nu(x).
           \label{eq:stable-polynomial}
\end{align}
Here a vertical strip has at most one box in each row, and we allow
the empty strip. The coefficient $\chi^\nu(\alpha)$ comes from the
character table of the fixed group $S_{|\nu|}$. Also,
$\binom{x}{a}=(x)_a/a!$ is a polynomial. Since $\nu\subseteq\mu$,
only $x_1,\ldots,x_k$ occur. Thus we never need the character table
of the growing group $S_n$ to compute $q_\mu$.

Formula \eqref{eq:stable-polynomial} is the classical stable
Specht-character formula of Macdonald
\cite[Example~I.7.14(b)]{Macdonald}; see also Narayanan, Paul,
Prasad, and Srivastava \cite[Section~3.1, equation~(15)]{NPPS}.
We now derive it using induction and Pieri's rule, so that we can
see why this particular polynomial gives the character.

For $\mu\vdash k$, consider
$M_n(\mu)=\operatorname{Ind}_{S_{n-k}\times S_k}^{S_n}
(\mathbf1\boxtimes S^\mu)$. Its underlying vector-space summands
are indexed by $k$-subsets. A subset contributes to the trace only
if it is invariant under the permutation, which means that it is
a union of cycles. To select a subset with cycle type
$\alpha\vdash k$, we choose $m_j(\alpha)$ of the $j$-cycles for
each $j$. There are $\prod_j\binom{C_j(\sigma)}{m_j(\alpha)}$
such choices, each contributing $\chi^\mu(\alpha)$. Summing over
the possible cycle types gives
\begin{equation}
 \chi_{M_n(\mu)}(\sigma)=F_\mu(C(\sigma)),\qquad
 \dim M_n(\mu)=\binom nk\dim S^\mu.
 \label{eq:induced-character-compact}
\end{equation}

Recall that a horizontal strip has at most one box in each column,
and its removal must leave a partition. Pieri's rule gives the
following decomposition when $n\ge2k$:
\begin{equation}
 M_n(\mu)\cong
 \bigoplus_{\substack{\nu\subseteq\mu\\\mu/\nu\text{ horizontal}}}
 S^{\nu[n]}.
 \label{eq:pieri-tail}
\end{equation}
To see the tail indexing, the condition that $\lambda/\mu$ is
horizontal is the interlacing condition
$\lambda_1\ge\mu_1\ge\lambda_2\ge\mu_2\ge\cdots$.
Set $\nu_i=\lambda_{i+1}$, so $\nu$ is the part below the first
row. Conversely, $n\ge2k$ ensures that every tail in the displayed
sum can be padded to a partition. Notice that the term $\nu=\mu$
is present, and every other term has a smaller tail. We can
therefore solve for the desired character recursively.

There is also a direct formula for this inversion. For the complete
and elementary symmetric functions, let $H(t)=\sum h_rt^r$ and
$E(t)=\sum e_rt^r$. Recall that $H(t)E(-t)=1$. By Pieri's rules,
multiplying $\sum_\nu u(\nu)s_\nu t^{|\nu|}$ by $H(t)$ sums over
horizontal strips, while multiplication by $E(-t)$ reverses that
operation by summing over vertical strips with alternating signs.
Comparing coefficients gives
\begin{equation}
 v(\mu)=\sum_{\mu/\nu\text{ horizontal}}u(\nu)
 \quad\Longrightarrow\quad
 u(\mu)=\sum_{\mu/\nu\text{ vertical}}
                 (-1)^{|\mu|-|\nu|}v(\nu).
 \label{eq:strip-inversion}
\end{equation}
Only degrees through $k$ are needed, so the series can be truncated
there.

\begin{proposition}\label{prop:stable-character}
For $\mu\vdash k$, $n\ge2k$, and $\sigma\in S_n$,
\begin{align}
 \chi^{\mu[n]}(\sigma)&=q_\mu(C_1(\sigma),\ldots,C_k(\sigma)),
 \label{eq:exact-stable-character}\\
 \dim S^{\mu[n]}&=\frac{\dim S^\mu}{k!}n^k+O_\mu(n^{k-1}).
 \label{eq:fixed-tail-dimension}
\end{align}
Equivalently, $q_\varnothing=1$ and
$q_\mu=F_\mu-\sum_{\nu\subsetneq\mu,\ \mu/\nu\text{ horizontal}}q_\nu$.
\end{proposition}
\begin{proof}
Take characters in \eqref{eq:pieri-tail} and apply the inversion
formula \eqref{eq:strip-inversion}. The resulting alternating sum is
exactly \eqref{eq:stable-polynomial}, so it gives the asserted
character formula.

For the dimension, evaluate at the identity permutation. Its cycle
counts are $(n,0,\ldots)$, and
$F_\nu(n,0,\ldots)=\dim S^\nu\binom n{|\nu|}$. This is a polynomial
in $n$ of degree $|\nu|$. In the sum defining $q_\mu$, the only term
of degree $k$ is the term with $\nu=\mu$; all the others have degree
at most $k-1$. Since the leading coefficient of $\binom nk$ is $1/k!$,
we obtain \eqref{eq:fixed-tail-dimension}.
\end{proof}
The classical formula in fact holds throughout the valid padded
range $n\ge k+\mu_1$ \cite[Section~2.1]{NPPS}. The simpler
sufficient range $n\ge2k$ used above is enough for our purposes.

For a transposed diagram, we use the standard identity
\begin{equation}
 \chi^{\lambda'}(\sigma)=\sgn(\sigma)\chi^\lambda(\sigma).
 \label{eq:transpose-character}
\end{equation}
See Macdonald \cite[Example~I.7.2]{Macdonald}. Here it is important
to distinguish $(\mu[n])'$ from $(\mu')[n]$: in the first expression
we transpose the entire padded diagram, not just the fixed tail.
Deleting the first column from the original diagram leaves a fixed
shape exactly when deleting the first row from its transpose does.

Before taking limits, let us ask whether a fixed cycle-polynomial
formula forces the diagram to have one of these shapes. We must
keep the sign factor in \eqref{eq:transpose-character}, since the
first finitely many cycle counts do not determine it. For example,
fix $r$ and take $n\ge2r+2$. An $n$-cycle and a permutation of cycle
type $(r+1,n-r-1)$ both have $C_1=\cdots=C_r=0$, but their signs
are opposite. Thus the row and column cases require different
character formulas.

We will use the Murnaghan--Nakayama rule. Recall that a
\emph{border strip}, also called a rim hook, is a skew Young diagram
connected through shared edges and containing no $2\times2$ block.
Its height $\operatorname{ht}$ is the number of rows it occupies
minus one. A removable strip must leave a Young diagram.

Suppose $\lambda\vdash n+L$, $\sigma\in S_n$, and $\gamma$ is an
$L$-cycle on $L$ new labels. The rule says that
\begin{equation}
 \chi^\lambda(\sigma\gamma)
 =\sum_{\substack{\nu\subseteq\lambda,\ |\nu|=n\\
                  \lambda/\nu\text{ a border strip}}}
       (-1)^{\operatorname{ht}(\lambda/\nu)}\chi^\nu(\sigma).
 \label{eq:mn-long-cycle}
\end{equation}
In other words, we remove a border strip with the same length as
$\gamma$, assign the sign determined by its height, and evaluate
the remaining character at $\sigma$. See Macdonald
\cite[Example~I.7.5]{Macdonald} or Kochhar--Wildon
\cite[Theorem~1.1]{KochharWildon}. In the proof below, we use this
as an identity for every $\sigma\in S_n$, not merely a calculation
at one permutation.

\begin{proposition}[Exact cycle polynomials and stabilized tails]
\label{prop:exact-tail-equivalence}
Let $\lambda_n\vdash n$ and fix $e\in\{0,1\}$. There exist a fixed integer
$r\ge1$, a fixed polynomial $P\in\R[x_1,\ldots,x_r]$, and an integer $N$
such that
\begin{equation}
 \chi^{\lambda_n}(\sigma)
 =\sgn(\sigma)^e P(C_1(\sigma),\ldots,C_r(\sigma))
 \qquad(n\ge N,\ \sigma\in S_n)
 \label{eq:exact-small-cycle-formula}
\end{equation}
if and only if there is a fixed partition $\mu$, possibly empty, such that
for every sufficiently large $n$,
\begin{equation}
 \lambda_n=
 \begin{cases}
  \mu[n],&e=0,\\
  (\mu[n])',&e=1.
 \end{cases}
 \label{eq:exact-tail-stabilization}
\end{equation}
Thus the diagram stabilizes after deletion of the first row when $e=0$,
and after deletion of the first column when $e=1$. The assertion is about
an exact identity for every permutation, with $P$ and $e$ independent of
$n$; no dimension-growth hypothesis is required.
\end{proposition}

\begin{proof}
Suppose first that \eqref{eq:exact-tail-stabilization} holds. Take
$P=q_\mu$ and $r\ge\max\{1,|\mu|\}$. The fixed-tail formula in
Proposition~\ref{prop:stable-character}, together with
\eqref{eq:transpose-character}, gives
\eqref{eq:exact-small-cycle-formula}. It remains to prove the
converse. The idea is to append a cycle too long to affect any of the
cycle counts appearing in $P$, then use Murnaghan--Nakayama to see what
this requires of the diagrams.

Fix $n\ge N$ and $L>r$, and let $\gamma$ be an $L$-cycle on new
labels, disjoint from $\sigma\in S_n$. Appending $\gamma$ leaves
$C_1,\ldots,C_r$ unchanged. It multiplies the sign by $(-1)^{L-1}$,
but the polynomial $P$ is the same at degrees $n$ and $n+L$. Hence
\begin{equation}
 \chi^{\lambda_{n+L}}(\sigma\gamma)
 =(-1)^{e(L-1)}\chi^{\lambda_n}(\sigma)
 \qquad(\sigma\in S_n).
 \label{eq:mn-stable-extension}
\end{equation}
Now expand the left side using \eqref{eq:mn-long-cycle}. Each remaining
partition $\nu$ occurs once, with coefficient $+1$ or $-1$. The
identity holds for every $\sigma\in S_n$, and the irreducible
characters are linearly independent. Thus the coefficient of
$\chi^{\lambda_n}$ must be $(-1)^{e(L-1)}$, and every other
coefficient must be zero. There can therefore be exactly one term,
with remaining diagram $\lambda_n$. In particular,
\begin{equation}
 \begin{gathered}
  \lambda_n\subseteq\lambda_{n+L},\qquad
  \lambda_{n+L}/\lambda_n\text{ is a border strip},\\
  \operatorname{ht}(\lambda_{n+L}/\lambda_n)
       \equiv e(L-1)\pmod2
  \qquad(n\ge N,\ L>r).
 \end{gathered}
 \label{eq:mn-nested-strips}
\end{equation}
Notice why it matters that we compared class functions. Character
values at one permutation may cancel, but coefficients of distinct
irreducible characters in a class-function identity cannot cancel.

We next use the geometry of these strips. Write $D_n$ for the diagram
of $\lambda_n$, and let $w=\lambda_{N,1}$ and $h=\lambda'_{N,1}$.
Increase $N$ to at least one if necessary. For every $m>N+r$,
$D_m/D_N$ is a border strip, so it cannot contain a $2\times2$ block.
This gives the bounds
\begin{equation}
 \lambda_{m,2}\le w+1,\qquad \lambda'_{m,2}\le h+1.
 \label{eq:mn-interior-bounds}
\end{equation}
Here we interpret a missing part as zero. To prove the first bound,
suppose the second row of $D_m$ extended through column $w+2$. Since
rows of a Young diagram are decreasing in length, all four boxes in
rows $1,2$ and columns $w+1,w+2$ would belong to $D_m$. None belongs
to $D_N$, whose width is $w$. This would put a $2\times2$ block in
$D_m/D_N$, a contradiction. Transposing the argument gives the second
bound. Consequently, every box outside both the first row and the
first column lies in the fixed finite rectangle
$\{2,\ldots,h+1\}\times\{2,\ldots,w+1\}$.

We claim that the boxes in this rectangle eventually stop changing.
By \eqref{eq:mn-nested-strips}, $D_a\subseteq D_b$ whenever $a\ge N$
and $b-a>r$. Thus, once a box occurs at an index $a\ge N$, it occurs
at every index $b>a+r$. For each position in the rectangle, either no
box ever occurs for $m>N+r$, or a box is present at every sufficiently
large index. There are only finitely many positions, so we can choose
one $M>N+r$ beyond all of their stabilization indices. For $m\ge M$,
only the first row and the first column can still change.

Now take $a,b\ge M$ such that $b-a>r$ is even. The skew diagram
$D_b/D_a$ is connected, and all of its boxes lie in the first row or
the first column. It cannot use both. Indeed, the new segment at the
right end of the first row and the new segment at the bottom of the
first column are disconnected; a path between them would have to pass
through boxes already in $D_a$. Hence the strip is entirely horizontal
or entirely vertical.

The choice of an even length now distinguishes these possibilities.
A horizontal strip has height zero and sign $+1$. A vertical strip of
length $b-a$ has height $b-a-1$ and sign $-1$. The congruence in
\eqref{eq:mn-nested-strips} therefore forces the new boxes to lie in
the first row when $e=0$, and in the first column when $e=1$.

So far we have compared indices separated by an even number. To
finish, we must include every sufficiently large index. Fix $a\ge M$.
For any $n>a+r$, choose $b>n+r$ with $b-a$ even. By
\eqref{eq:mn-nested-strips}, $D_a\subseteq D_n\subseteq D_b$.
If $e=0$, the outer two diagrams have exactly the same boxes below
the first row. The intermediate diagram must have those same boxes.
Call their fixed shape $\mu$. Since $D_n$ has $n$ boxes in total,
$\lambda_n=\mu[n]$ for every $n>a+r$.

If $e=1$, the same inclusions fix all boxes outside the first column.
After transposition, these form one fixed first-row tail $\mu$, giving
$\lambda_n=(\mu[n])'$. Thus stabilization holds along the whole
sequence, not just along either parity of $n$.
\end{proof}

For direct sums, we will allow a sum of an unsigned polynomial and
a polynomial multiplied by the sign. We should therefore also check
what this more general expression means for an irreducible
character. The next corollary shows that it still gives only the
row or the column alternative.

\begin{corollary}
\label{cor:row-column-polynomials}
Let $\lambda_n\vdash n$. There exist fixed polynomials $f,g$ in a fixed
finite collection of cycle-count variables such that, for every
sufficiently large $n$ and every $\sigma\in S_n$,
\begin{equation}
 \chi^{\lambda_n}(\sigma)
 =f(C(\sigma))+\sgn(\sigma)g(C(\sigma))
 \label{eq:exact-signed-cycle-polynomials}
\end{equation}
if and only if, for one fixed partition $\mu$, either $\lambda_n=\mu[n]$
eventually or $\lambda_n=(\mu[n])'$ eventually. In the first case
$(f,g)=(q_\mu,0)$, and in the second $(f,g)=(0,q_\mu)$, as identities in
the common polynomial ring. In particular, the orientation cannot
alternate indefinitely.
\end{corollary}

\begin{proof}
Give the variable $x_j$ weight $j$, and choose $d$ at least as large
as the weighted degrees of $f$ and $g$. We use the fact that the
polynomials $q_\mu$ with $|\mu|\le d$ form a basis for the
polynomials of weighted degree at most $d$; see
\cite[Corollary~3.2]{NPPS}. Let us recall how this follows from the
formulas already given.

For $|\nu|=k$, the coefficients of $F_\nu$ in the binomial basis in
\eqref{eq:induced-polynomial} are the entries of the character table
of $S_k$. That table is invertible, so the highest-weight parts of
these $F_\nu$ form a basis at weight $k$. In
\eqref{eq:stable-polynomial}, $q_\mu$ differs from $F_\mu$ only by
terms of strictly smaller weight. We can therefore express the
highest-weight part of any polynomial using the $q_\mu$, subtract
it, and repeat at the next lower weight. This proves the basis
assertion by descending induction.

Write $f=\sum_{|\mu|\le d}a_\mu q_\mu$ and
$g=\sum_{|\mu|\le d}b_\mu q_\mu$. These coefficients do not depend
on $n$. By Proposition~\ref{prop:stable-character} and
\eqref{eq:transpose-character}, equation
\eqref{eq:exact-signed-cycle-polynomials} becomes an expansion in
$\chi^{\mu[n]}$ and $\chi^{(\mu[n])'}$.

For $n>2d+1$, all the irreducible characters in this expansion are
distinct. To see why the two families do not overlap, a row-padded
diagram has first row at least $n-d>d+1$, whereas a column-padded
diagram has first row at most $d+1$. Within either family, different
tails give different diagrams. Character orthogonality now forces
exactly one coefficient to be $1$ and all the others to be zero,
since the character on the left is irreducible. The coefficients
are fixed, so the same tail and the same orientation occur for every
sufficiently large $n$. This also identifies $f$ and $g$ as claimed.
The reverse implication follows from
Proposition~\ref{prop:exact-tail-equivalence}.
\end{proof}

\begin{remark}\label{rem:exact-identity}
Notice that the necessity proof of
Proposition~\ref{prop:exact-tail-equivalence} never used the degree
or the coefficients of $P$. It only used that one fixed function of
the small cycle counts is unchanged when a long cycle is appended.
Thus the same argument works with an arbitrary fixed function in
place of $P$, provided the sign factor is still fixed. Once the
diagram has stabilized, the fixed-tail formula gives a polynomial
description.

On the other hand, we did use an \emph{exact} identity at every
permutation. Convergence in distribution, or an identity with an
$o_{\Pp}(1)$ error, does not give that identity. The converse
therefore does not make the small-hook condition below necessary,
and it does not impose a diagram restriction on the centered
Gaussian theorem.
\end{remark}

We now turn from exact identities to limiting laws. The fixed-tail
formula tells us to examine the small cycle counts, while the
transpose identity tells us to keep the parity sign as well. The
following lemma gives their joint limit.

\begin{lemma}\label{lem:cycles-parity}
For every fixed $K$,
\begin{equation}
 (C_1^{(n)},\ldots,C_K^{(n)},\varepsilon_n)
 \dto(P_1,\ldots,P_K,\varepsilon),
 \label{eq:cycles-parity}
\end{equation}
where the variables on the right are independent, $P_j\sim\Pois(1/j)$,
and $\Pp(\varepsilon=1)=\Pp(\varepsilon=-1)=1/2$. Every fixed joint
polynomial moment also converges.
\end{lemma}

\begin{proof}
We begin with factorial moments. Fix nonnegative integers $m_j$ and
let $L=\sum_{j=1}^Kjm_j$. The product of falling factorials counts
ordered selections of $m_j$ distinct $j$-cycles for each $j$. These
cycles use $L$ labels. Count the selections on disjoint labels, then
complete the permutation on the remaining labels. This gives
\begin{equation}
 \Ee\prod_{j=1}^K(C_j^{(n)})_{m_j}=\prod_{j=1}^Kj^{-m_j}
 \qquad(n\ge L).
 \label{eq:cycle-factorial-moments}
\end{equation}
The right side is exactly the corresponding factorial moment of
independent Poisson variables with means $1/j$.

To include the sign, insert $\varepsilon_n$ in the same count. When
$n\ge L+2$, at least two labels remain after the selected cycles have
been fixed. A transposition of two remaining labels pairs each
completion with one of the opposite sign. Their contributions cancel,
so
\begin{equation}
 \Ee\!\left[\varepsilon_n\prod_{j=1}^K(C_j^{(n)})_{m_j}\right]=0
 \qquad(n\ge L+2).
 \label{eq:signed-factorial-moments}
\end{equation}
Since $\varepsilon_n^2=1$, these two identities determine every joint
polynomial moment, and each is eventually equal to its proposed
limiting value. A linear combination of the proposed limiting
variables has moment generating function
$\exp(\sum_j(e^{vt_j}-1)/j)\cosh(vu)$, which is finite near zero.
Its moments therefore determine its law. The method of moments,
followed by Cram\'er--Wold, gives the asserted joint convergence.
For quantitative approximations of the small cycle counts, see
\cite{ArratiaTavare}.
\end{proof}

The limiting sign is independent of the small cycle counts, rather
than a function of them. Combining this fact with
\eqref{eq:transpose-character} will give the limiting character for
a transposed fixed tail.

The general square-character bound already shows that the normalized
character at the square tends to zero whenever the ordinary character
converges. For a fixed tail we can say more: we can also identify the
limit of the unnormalized character at the square.

Recall what happens when a cycle is squared. An odd cycle remains a
cycle of the same length, while an even cycle splits into two cycles
of half the length. Hence
\begin{equation}
 C_j(\sigma^2)=
 \begin{cases}C_j(\sigma)+2C_{2j}(\sigma),&j\text{ odd},\\
               2C_{2j}(\sigma),&j\text{ even}.
 \end{cases}
 \label{eq:square-cycles}
\end{equation}
Thus the cycle counts of the square are functions of the original
cycle counts. They are not a new cycle vector from an independent
uniform permutation.

\begin{proposition}
\label{prop:tail-limits}
For fixed $\mu\vdash k$,
\[
 \chi^{\mu[n]}(\sigma_n)\dto q_\mu(P),\qquad
 \chi^{(\mu[n])'}(\sigma_n)\dto\varepsilon q_\mu(P).
\]
These limits hold jointly for any finite collection of fixed tails. Both corresponding character sequences
at $\sigma_n^2$ are tight and have the same limit
\[
 q_\mu(Y_1,\ldots,Y_k),\qquad
 Y_j=\begin{cases}P_j+2P_{2j},&j\text{ odd},\\2P_{2j},&j\text{ even},\end{cases}
\]
where the $P_j$ through $2k$ are independent with means $1/j$.
\end{proposition}

\begin{proof}
By Proposition~\ref{prop:stable-character}, each fixed-tail character
is the evaluation of one fixed polynomial. Apply
Lemma~\ref{lem:cycles-parity} and continuity of polynomial evaluation
to obtain its limit. For a transposed diagram,
\eqref{eq:transpose-character} multiplies the same polynomial by the
same parity sign. Using a single cycle vector large enough for all
the tails gives the joint assertion.

Now consider $\sigma_n^2$. Its sign is always $+1$, so a diagram and
its transpose have the same character at this permutation. Substitute
\eqref{eq:square-cycles} into the fixed-tail formula, and apply the
cycle-count limit through $2k$. This gives the stated polynomial in
the variables $Y_j$. Since the resulting characters converge in
distribution, they are also tight.
\end{proof}

We next consider diagrams with small hooks. Fixed character
polynomials do not describe these growing diagrams, so we use a
different argument. We will find one event whose probability tends
to one and on which every constituent character is zero, both at
the permutation and at its square.

Return to the Murnaghan--Nakayama rule
\eqref{eq:mn-long-cycle}. If there is no removable border strip of
the required cycle length, its sum is empty and the character is
zero. Since we may consider the cycles in any order, we can start
with a longest cycle.

\begin{lemma}\label{lem:long-cycle}
If $\lambda\vdash n$ and a permutation $\sigma$ has a cycle longer than
$h(\lambda)$, then $\chi^\lambda(\sigma)=0$.
\end{lemma}

\begin{proof}
Let us first bound the size of a border strip inside $\lambda$.
Follow its boxes from the southwest end to the northeast end. Each
step goes upward or to the right, so a strip occupying $r$ rows and
$c$ columns has $r+c-1$ boxes. Since it lies inside $\lambda$, this
is at most $\ell(\lambda)+\lambda_1-1=h(\lambda)$.

Now apply Murnaghan--Nakayama to the cycle in the hypothesis. Its
length exceeds the size of every possible border strip, so there is
no strip of the required length to remove. The sum in the rule is
empty, and hence $\chi^\lambda(\sigma)=0$.
\end{proof}

Let $L_n$ be the longest cycle length of a uniform permutation in
$S_n$, and write $p_{n,b}=\Pp(L_n\le b)$. For an integer
$0\le b\le n$, we have
\begin{equation}
 p_{n,b}\le\frac bn.
 \label{eq:longest-cycle-bound}
\end{equation}
To prove this, notice that on $\{L_n\le b\}$ all $n$ labels belong
to cycles of length at most $b$. Thus
$n\mathbf1_{\{L_n\le b\}}\le\sum_{j=1}^b jC_j^{(n)}$.
The expectation of the sum on the right is $b$, so taking
expectations gives the claimed bound. We set $p_{n,b}=1$ for
$b\ge n$.

\begin{proposition}
\label{prop:small-hook}
Let $R_n$ be any real representation of $S_n$ whose constituent diagrams
have largest hook at most $H_n$, and let $\eta_n$ be its character. Then
\begin{equation}
 \Pp\bigl(\eta_n(\sigma_n)\ne0\ \text{or }
                 \eta_n(\sigma_n^2)\ne0\bigr)
 \le p_{n,2H_n}\le\min\{1,2H_n/n\}.
 \label{eq:simultaneous-vanishing}
\end{equation}
In particular, if $H_n=o(n)$, both characters converge to zero in
probability, without any restriction on the constituent multiplicities.
For a single irreducible diagram with $h(\lambda_n)=o(n)$, its dimension
also tends to infinity.
\end{proposition}

\begin{proof}
Consider the event $\{L_n>2H_n\}$. The original permutation has a
cycle longer than $H_n$. Squaring splits this cycle into at most two
cycles, so the square has a cycle of length at least $L_n/2>H_n$ as
well. Every constituent of $R_n$ has largest hook at most $H_n$.
Lemma~\ref{lem:long-cycle} therefore makes every constituent
character zero at both permutations on this same event.

The character of $R_n$ is the sum of these constituent characters,
with their multiplicities. A sum of zero terms is still zero, no
matter what those multiplicities are. Thus either character can be
nonzero only on the complementary event $\{L_n\le2H_n\}$, proving
the first inequality. The second follows from
\eqref{eq:longest-cycle-bound}. If $R_n=0$, both characters are
identically zero and we take $H_n=0$.

For a single irreducible diagram, the hook-length formula gives
\[
 \dim S^{\lambda_n}=\frac{n!}{\prod_{u\in\lambda_n}h_u}
       \ge\frac{n!}{h(\lambda_n)^n}
       \ge\left(\frac{n}{e\,h(\lambda_n)}\right)^n\longrightarrow\infty.
\]
In the first inequality, we used that each box hook is at most the
upper-left hook. In the second, we used $n!\ge(n/e)^n$, which follows
by integrating $\log x$ from $1$ to $n$. Finally,
$h(\lambda_n)=o(n)$ makes $n/(e\,h(\lambda_n))$ tend to infinity.
This proves the dimension claim.
\end{proof}

Notice that we have not merely bounded the character after dividing
by the dimension. The ordinary character itself is exactly zero on
an event whose probability tends to one. Furthermore, the event
works for all constituents at once. This is why arbitrary
multiplicities cause no difficulty; we do not have to add separate
probability bounds for the individual summands.

\begin{figure}[H]
\centering
\begin{tikzpicture}[font=\small]
\node[align=center] at (2.1,.8) {Fixed first-row tail\\$\mu=(2,1)$};
\node[align=center] at (7.0,.8) {Transposed fixed tail\\$\mu[n]\,{}^{\prime}$};
\node[align=center] at (12.1,.8) {Sublinear width and height\\$h(\lambda_n)=o(n)$};
\begin{scope}[shift={(0.65,0)},x=0.48cm,y=0.48cm]
\draw[line width=.4pt,fill=white] (0,0) rectangle (1,-1);
\draw[line width=.4pt,fill=white] (1,0) rectangle (2,-1);
\draw[line width=.4pt,fill=white] (2,0) rectangle (3,-1);
\draw[line width=.4pt,fill=white] (3,0) rectangle (4,-1);
\draw[line width=.4pt,fill=white] (4,0) rectangle (5,-1);
\draw[line width=.4pt,fill=white] (5,0) rectangle (6,-1);
\draw[line width=.4pt,fill=black!18] (0,-1) rectangle (1,-2);
\draw[line width=.4pt,fill=black!18] (1,-1) rectangle (2,-2);
\draw[line width=.4pt,fill=black!18] (0,-2) rectangle (1,-3);
\end{scope}
\begin{scope}[shift={(6.2,0)},x=0.43cm,y=0.43cm]
\draw[line width=.4pt,fill=white] (0,0) rectangle (1,-1);
\draw[line width=.4pt,fill=black!18] (1,0) rectangle (2,-1);
\draw[line width=.4pt,fill=black!18] (2,0) rectangle (3,-1);
\draw[line width=.4pt,fill=white] (0,-1) rectangle (1,-2);
\draw[line width=.4pt,fill=black!18] (1,-1) rectangle (2,-2);
\draw[line width=.4pt,fill=white] (0,-2) rectangle (1,-3);
\draw[line width=.4pt,fill=white] (0,-3) rectangle (1,-4);
\draw[line width=.4pt,fill=white] (0,-4) rectangle (1,-5);
\draw[line width=.4pt,fill=white] (0,-5) rectangle (1,-6);
\end{scope}
\begin{scope}[shift={(11.2,0)},x=0.55cm,y=0.55cm]
\draw[line width=.4pt,fill=black!18] (0,0) rectangle (1,-1);
\draw[line width=.4pt,fill=black!18] (1,0) rectangle (2,-1);
\draw[line width=.4pt,fill=black!18] (2,0) rectangle (3,-1);
\draw[line width=.4pt,fill=black!18] (0,-1) rectangle (1,-2);
\draw[line width=.4pt,fill=white] (1,-1) rectangle (2,-2);
\draw[line width=.4pt,fill=white] (2,-1) rectangle (3,-2);
\draw[line width=.4pt,fill=black!18] (0,-2) rectangle (1,-3);
\draw[line width=.4pt,fill=white] (1,-2) rectangle (2,-3);
\draw[line width=.4pt,fill=white] (2,-2) rectangle (3,-3);
\end{scope}

\draw[-{Stealth},line width=.7pt] (3.7,-.24)--(4.55,-.24);
\node[font=\scriptsize,align=center] at (4.12,-.83) {row grows};
\draw[-{Stealth},line width=.7pt] (6.415,-2.72)--(6.415,-3.25);
\node[font=\scriptsize,anchor=west] at (6.72,-3.0) {column grows};
\draw[<->,line width=.6pt] (11.2,-1.94)--(12.85,-1.94);
\node[font=\scriptsize] at (12.025,-2.22) {width};
\draw[<->,line width=.6pt] (13.1,0)--(13.1,-1.65);
\node[font=\scriptsize,rotate=90] at (13.4,-.82) {height};
\node[align=center] at (2.3,-3.95) {$q_\mu(P)$};
\node[align=center] at (7,-3.95) {$\varepsilon q_\mu(P)$};
\node[align=center] at (12.1,-3.95) {$0$};
\end{tikzpicture}

\caption{Three geometric inputs and their character limits. All pictures have nine boxes; they show finite members, not asymptotic thresholds. In the left and middle pictures the three shaded boxes stay fixed while the long row or column grows. In the square, the shaded upper-left hook has length $3+3-1=5$; along growing squares its length is negligible compared with the area. Transposition of the whole padded diagram introduces the sign $\varepsilon$.}\label{fig:regimes}
\end{figure}
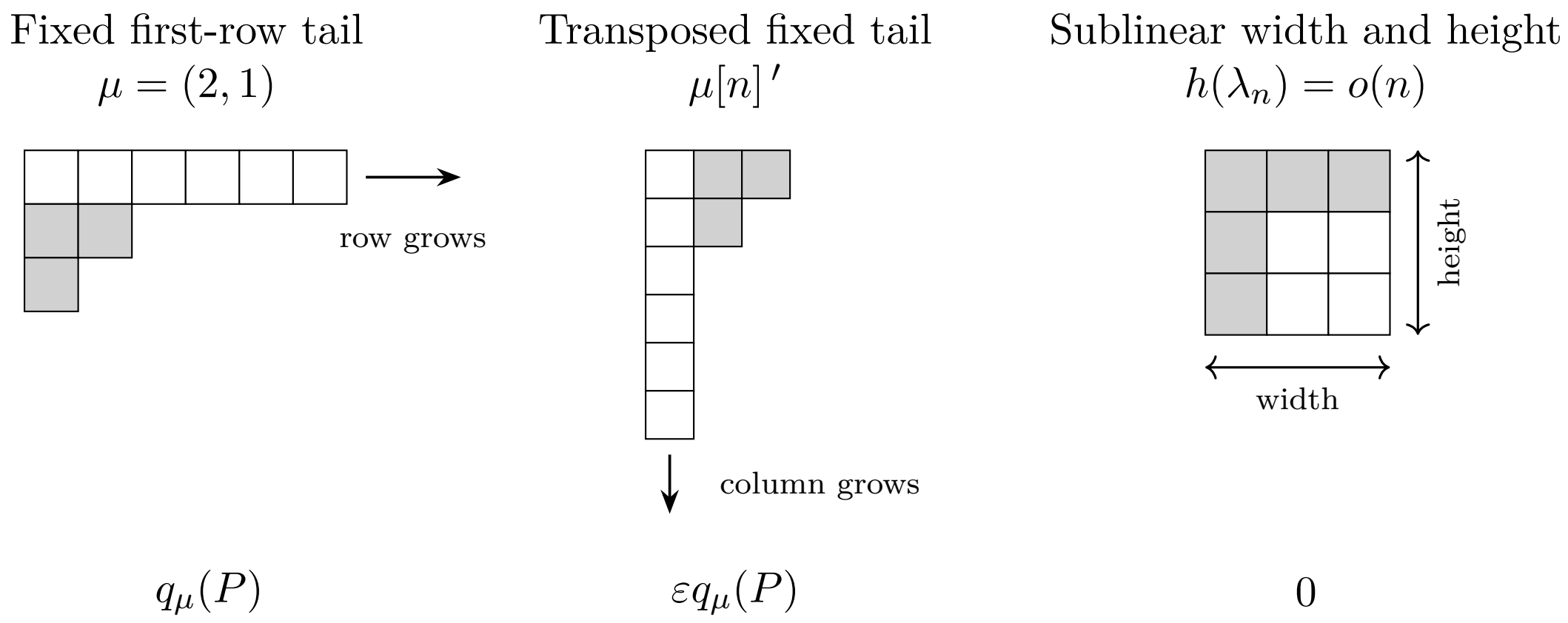

We can now apply Theorem~\ref{thm:character-equivalence} to these
character laws. Figure~\ref{fig:regimes} shows the three types of
diagrams. In the next theorem, the fixed tails and their fixed
multiplicities determine the limiting character, while the
small-hook remainder may have arbitrary multiplicities.

\begin{theorem}[Explicit character-plus-Gaussian limits]\label{thm:diagram}
Suppose $d_n\to\infty$. Let $\mathcal T$ be a finite set of fixed partitions,
possibly including $\varnothing$, and let $a_\mu,b_\mu$ be fixed nonnegative
integers. Assume that, for all sufficiently large $n$,
\begin{equation}
 \rho_n\cong
 \bigoplus_{\mu\in\mathcal T}(S^{\mu[n]})^{\oplus a_\mu}
 \ \oplus\!
 \bigoplus_{\mu\in\mathcal T}(S^{(\mu[n])'})^{\oplus b_\mu}
 \ \oplus R_n,
 \qquad \frac{H_n}{n}\longrightarrow0,
 \tag{D}\label{eq:diagram-condition}
\end{equation}
where $H_n=\max\{h(\lambda):S^\lambda\text{ occurs in }R_n\}$, with
$H_n=0$ when $R_n=0$. There is no restriction on the multiplicities or the
number of constituents of $R_n$, apart from the displayed maximum.

Take $K\ge1$ at least as large as every $|\mu|$ in $\mathcal T$, and let
$P_1,\ldots,P_K,\varepsilon,G$ be independent, with
$P_j\sim\Pois(1/j)$, $\Pp(\varepsilon=1)=\Pp(\varepsilon=-1)=1/2$, and
$G\sim N(0,1)$. Define
\begin{equation}
 \Xi=\sum_{\mu\in\mathcal T}a_\mu q_\mu(P)
          +\varepsilon\sum_{\mu\in\mathcal T}b_\mu q_\mu(P),
 \qquad P=(P_1,\ldots,P_K).
 \label{eq:explicit-Xi}
\end{equation}
Then $\chi_n(\sigma_n)\dto\Xi$ and
$\chi_n(\sigma_n^2)/d_n\to0$ in $L^2$. In fact, the unnormalized
square characters $\chi_n(\sigma_n^2)$ are tight.

For tests satisfying~\eqref{eq:test-normalization}, write
$T_n=\Tr(A_nO_nQ_nO_n^T)$. Then
\begin{equation}
 \bigl(C_1^{(n)},\ldots,C_K^{(n)},\varepsilon_n,
       \chi_n(\sigma_n),T_n-s_n\chi_n(\sigma_n)\bigr)
 \dto
 \bigl(P_1,\ldots,P_K,\varepsilon,\Xi,\sqrt{1-s^2}\,G\bigr).
 \label{eq:diagram-joint}
\end{equation}
In particular,
\begin{equation}
 T_n\dto s\Xi+\sqrt{1-s^2}\,G.
 \label{eq:diagram-limit}
\end{equation}
\end{theorem}

\begin{proof}
Let $\eta_n$ be the character of $R_n$. Characters add under direct
sums, so \textup{(D)} gives
$\chi_n(\sigma_n)=f(C^{(n)})+\varepsilon_ng(C^{(n)})+\eta_n(\sigma_n)$,
where $f=\sum a_\mu q_\mu$ and $g=\sum b_\mu q_\mu$.
The polynomials are fixed because the tails and their multiplicities
are fixed. Lemma~\ref{lem:cycles-parity} gives the joint limit of the
first two terms with the small cycle counts and parity. Meanwhile,
Proposition~\ref{prop:small-hook} gives
$\eta_n(\sigma_n)\to0$ in probability. Hence
$\chi_n(\sigma_n)\dto\Xi$, jointly with
$(C_1^{(n)},\ldots,C_K^{(n)},\varepsilon_n)$.

Now apply Theorem~\ref{thm:character-equivalence}, taking this
small-cycle and parity vector as $W_n$. We obtain the normalized
square-character convergence in $L^2$, the joint limit
\eqref{eq:diagram-joint}, and the uncentered limit
\eqref{eq:diagram-limit}. No further trace condition is needed.

Lastly, we check the stronger assertion about the unnormalized
square character. The fixed-tail part is tight by
Proposition~\ref{prop:tail-limits}, and the remainder tends to zero
in probability by Proposition~\ref{prop:small-hook}. Their sum is
therefore tight, which completes the proof.
\end{proof}

The character limits in the theorem follow from the diagram
condition; they are not additional hypotheses. Notice also that
all summands in \eqref{eq:explicit-Xi} use the same Poisson variables
and the same sign. They come from one permutation, so they must not
be replaced by independent copies. The Gaussian variable, by
contrast, is independent of the whole collection. When $|s|=1$,
its coefficient is zero.

The fixed-multiplicity assumption is used to identify the ordinary
character limit. It is not needed for the centered conditional
central limit theorem: Theorem~\ref{thm:universality} allows
arbitrary other diagrams and multiplicities as long as
$t_n=o(d_n)$. To identify the unnormalized square-character law, we
use the Poisson coordinates through $2K$. Those extra coordinates
are not needed for the normalized square bound or for the joint
vector in \eqref{eq:diagram-joint}. We may also choose a larger
fixed $K$ from the start. The polynomials simply ignore the extra
coordinates, and Lemma~\ref{lem:cycles-parity} gives the same
limits after projection onto a shorter vector.

\begin{corollary}\label{cor:irreducible}
Let $\rho_n$ be an orthogonal realization of $S^{\lambda_n}$. Each of the
following hypotheses implies $d_n\to\infty$ and the limit
\eqref{eq:diagram-limit}, for tests satisfying~\eqref{eq:test-normalization}:
\[
\begin{array}{c|c}
 \text{Condition on }\lambda_n & \Xi\\ \hline
 \lambda_n=\mu[n],\quad\mu\ne\varnothing\text{ fixed}
       &q_\mu(P_1,\ldots,P_{|\mu|})\\[2pt]
 \lambda_n=(\mu[n])',\quad\mu\ne\varnothing\text{ fixed}
       &\varepsilon q_\mu(P_1,\ldots,P_{|\mu|})\\[2pt]
 \lambda_{n,1}+\lambda'_{n,1}=o(n)&0
\end{array}
\]
\end{corollary}

\begin{proof}
Suppose first that the diagram has a fixed nonempty tail, or is the
transpose of such a diagram. Then \textup{(D)} holds with $R_n=0$,
and \eqref{eq:fixed-tail-dimension} gives $d_n\to\infty$.
Theorem~\ref{thm:diagram} gives the corresponding character law and
randomized-basis limit.

In the small-hook case, take the entire representation to be $R_n$
and set every fixed-tail multiplicity to zero.
Proposition~\ref{prop:small-hook} gives both dimension growth and
$\Xi=0$. Applying the same theorem proves the remaining case.
\end{proof}

The corollary identifies $\Xi$ for three types of irreducible
diagrams. Recall that the entry limits already hold for every
growing-dimensional irreducible sequence by
Corollary~\ref{cor:all-irreducible-entries}. For the reducible
representations in \textup{(D)}, character convergence also allows
us to remove the centering in Corollary~\ref{thm:entries}. Thus the
geometric assumptions here determine the ordinary character law;
they are not what makes the centered entry fluctuations Gaussian.

For the centered conclusion, the fixed-tail multiplicities in
\textup{(D)} may even grow with $n$, provided the total trivial and
sign multiplicity is $o(d_n)$. This follows from
Theorem~\ref{thm:universality}, which does not require a fixed-tail
or small-hook assumption at all. The next example explains why
this does not automatically give an ordinary character limit or
allow us to omit diagonal centering.

\begin{example}[A growing multiple of the standard representation]
\label{ex:standard-multiple}
Consider $\rho_n=(S^{(n-1,1)})^{\oplus n}$, whose dimension is
$d_n=n(n-1)$. Its ordinary character is
$\chi_n(\sigma_n)=n(C_1^{(n)}-1)$, which is not tight. Nevertheless,
the normalized characters are
\[
 \frac{\chi_n(\sigma_n)}{d_n}=\frac{C_1^{(n)}-1}{n-1},\qquad
 \frac{\chi_n(\sigma_n^2)}{d_n}
       =\frac{C_1^{(n)}+2C_2^{(n)}-1}{n-1}
\]
and both tend to zero in probability. There are no trivial or sign
summands, so $t_n=0$ and Theorem~\ref{thm:universality} gives the
centered conditional limit.

At entry scale, however,
$\chi_n(\sigma_n)/\sqrt{d_n}\dto P_1-1$ rather than zero.
Thus, for fixed $r$, the joint form of Corollary~\ref{thm:entries}
gives
\[
 \left(\sqrt{d_n}\,(O_nQ_nO_n^T)_{jj}\right)_{j=1}^r
 \dto (P_1-1+G_j)_{j=1}^r.
\]
Here $G_1,\ldots,G_r$ are independent standard normals, independent
of $P_1$. Each diagonal entry has the same additional Poisson
shift. Subtracting $\chi_n(\sigma_n)/\sqrt{d_n}$ from every scaled
diagonal entry removes that common term and leaves independent
standard Gaussian limits.

The common Poisson term is exactly the conditional mean that
survives at entry scale. Thus it is consistent with the centered
theorem, which subtracts this term. The example also shows why
$s_n\to0$ alone is not enough to discard
$s_n\chi_n(\sigma_n)$ when the character is not tight.
\end{example}

\begin{example}[The regular representation]\label{ex:regular}
Take the regular representation. Its dimension
is $d_n=n!$, and its character is
$\chi_n(\sigma)=n!\mathbf1_{\{\sigma=e\}}$, since a nonidentity
left translation fixes no basis vector. Thus
$\chi_n(\sigma_n)\to0$ in probability, even though
$\Ee[\chi_n(\sigma_n)^2]=n!$. Theorem~\ref{thm:character-equivalence}
applies with $\Xi=0$, and Corollary~\ref{thm:entries} gives
uncentered Gaussian entries.

Taking the inner product with the regular character shows that
$S^\lambda$ occurs with multiplicity $f^\lambda$. In particular,
$t_n=2$. This representation does not satisfy \textup{(D)}: the
standard constituent has multiplicity $n-1$ and largest hook $n$.
It therefore fits neither the fixed-multiplicity part nor the
small-hook remainder. This shows that the multiplicity criterion
applies to more representations than the sufficient diagram
condition used to compute explicit character laws.
\end{example}

Let us finish with some explicit calculations and then discuss the
limits of the diagram assumptions. We begin with a fixed tail whose
padded diagram is neither a two-row diagram nor a hook.

Take $\mu=(2,1)$. The corresponding $S_3$ character has values
$2,0,-1$ on cycle types $1^3,(2,1),(3)$, respectively. Substituting
these into \eqref{eq:induced-polynomial} gives
$F_{(2,1)}=2\binom{x_1}{3}-x_3$. The diagrams left after removing a
vertical strip are $(2,1),(2),(1,1),(1)$. Hence
\[
 q_{(2,1)}=F_{(2,1)}-F_{(2)}-F_{(1,1)}+F_{(1)}
 =2\binom{x_1}{3}-2\binom{x_1}{2}+x_1-x_3.
\]
Notice that the $x_2$ terms cancel. To find the dimension, evaluate
at the identity, whose cycle vector is $(n,0,\ldots)$; this gives
$d_n=n(n-2)(n-4)/3$. As a numerical example, take
$(1)(2)(3)(4\,5)(6\,7\,8)\in S_8$. Its first three cycle counts
are $(3,1,1)$, so substitution gives character $-2$, while the
dimension is $64$.

We can now replace the cycle counts by their Poisson limits.
For $\lambda_n=(n-3,2,1)$, this gives
\begin{equation}
 \Xi=2\binom{P_1}{3}-2\binom{P_1}{2}+P_1-P_3
     =\frac{P_1(P_1-2)(P_1-4)}3-P_3.
 \label{eq:tail21-limit}
\end{equation}
The general theorem then supplies the randomized-basis limit;
there is no separate trace-limit condition to check.
Table~\ref{tab:small-tails} lists the other first fixed-tail
polynomials.

\begin{table}[H]
\centering\small
\caption{Fixed-tail laws from~\eqref{eq:stable-polynomial}. Write
$B_r=\binom{x_1}{r}$. Each row is valid for $n\ge2|\mu|$; the same
polynomial with an independent common sign applies to the transposed
padded diagram.}\label{tab:small-tails}
\renewcommand{\arraystretch}{1.55}
\begin{tabularx}{\linewidth}{@{}c c >{\raggedright\arraybackslash}X@{}}
\toprule
Tail $\mu$ & Padded diagram & $q_\mu(x)$\\
\midrule
$\varnothing$&$(n)$&$1$\\
$(1)$&$(n-1,1)$&$x_1-1$\\
$(2)$&$(n-2,2)$&$B_2-x_1+x_2$\\
$(1,1)$&$(n-2,1,1)$&$B_2-x_1+1-x_2$\\
$(3)$&$(n-3,3)$&$B_3-B_2+x_1x_2-x_2+x_3$\\
$(2,1)$&$(n-3,2,1)$&$2B_3-2B_2+x_1-x_3$\\
$(1,1,1)$&$(n-3,1,1,1)$&$B_3-B_2+x_1-1-x_1x_2+x_2+x_3$\\
$(2,2)$&$(n-4,2,2)$&$2B_4-2B_3+B_2+2\binom{x_2}{2}-x_2-(x_1-1)x_3$\\
\bottomrule
\end{tabularx}
\end{table}

There is a useful check on these polynomials. Character
orthogonality and \eqref{eq:cycle-factorial-moments} give
\begin{equation}
 \Ee[q_\mu(P)q_\nu(P)]=\mathbf1_{\{\mu=\nu\}}.
 \label{eq:poisson-orthogonality}
\end{equation}
Indeed, the expectation of a fixed polynomial in cycle counts is
eventually exactly its Poisson expectation. Apply this to
$q_\mu q_\nu$, then use character orthogonality at a sufficiently
large degree. Taking the other tail to be empty gives mean zero
for every nonempty tail, and taking the two tails equal gives
variance one. This does not imply moment convergence for a
small-hook remainder. Its second moment can remain one even
while its character tends to zero in probability.

For a transpose pair,
$S^{(n-1,1)}\oplus S^{(2,1^{n-2})}$, the limiting character is
$\Xi=(1+\varepsilon)(P_1-1)$. It is not a sum of independent
character copies: both summands use the same Poisson variable and
sign. Similarly, the defining permutation representation
$S^{(n)}\oplus S^{(n-1,1)}$ has $\Xi=P_1$. This recovers Tsou's
Poisson-plus-Gaussian law \cite{TsouRandomized}.

For a small-hook example, take the square diagram $(m^m)$. Here
$n=m^2$ and $H=2m-1$, so the probability that either character is
nonzero is at most $p_{m^2,4m-2}\le(4m-2)/m^2$.
We can compute the sufficient-event probability more closely by
conditioning on the length of the cycle containing label $1$.
That length is uniform on $\{1,\ldots,n\}$, giving the recurrence
\begin{equation}
 p_{0,b}=1,\qquad
 p_{n,b}=\frac1n\sum_{j=1}^{\min(b,n)}p_{n-j,b}.
 \label{eq:restricted-cycle-recurrence}
\end{equation}
For $m=10,20,30$, the elementary bounds are
$0.380000,0.195000,0.131111$, respectively. The recurrence gives
approximately $0.106160$, $2.670958\times10^{-4}$, and
$1.237311\times10^{-7}$. These still bound the probability through
a sufficient event; they are not the actual probabilities that
the characters are nonzero. The decimals are rounded evaluations,
not rigorous enclosures.

We should keep two different conclusions separate. Condition
\textup{(D)} is sufficient to identify the ordinary character law,
but it does not classify every sequence with a character limit.
The row-or-column classification in
Corollary~\ref{cor:row-column-polynomials} concerns exact identities
at every permutation. It is not a classification of convergence
under a uniform permutation. A small-hook remainder disappears
in probability, not as an exact class function. The centered
Gaussian conclusion, on the other hand, is characterized for all
real representations by Theorem~\ref{thm:universality}.

Even bounded tail size need not give a character limit. For example,
alternate the diagrams $(n-2,2)$ and $(n-2,1,1)$. Their limiting
polynomials have different laws: Table~\ref{tab:small-tails} gives
$q_{(2)}(P)\ge-1$, whereas $q_{(1,1)}(P)$ is unbounded below.
Thus the two subsequences have different limits, despite the
tightness in \eqref{eq:irreducible-tightness}. Nevertheless,
Corollary~\ref{cor:all-irreducible-entries} gives the centered
Gaussian and uncentered finite-entry limits along the whole
alternating sequence.

The calculations above do not identify the ordinary character laws
for growing hooks or for two-row diagrams with both rows
proportional to $n$. The general results already give their
normalized square-character convergence and Gaussian entry limits.
What remains is the ordinary character law. For hooks, the identity
$S^{(n-k,1^k)}\cong\Lambda^k S^{(n-1,1)}$ gives the exact formula
\begin{equation}
 \chi^{(n-k_n,1^{k_n})}(\sigma_n)
 =[z^{k_n}]\frac{1}{1+z}\prod_{j=1}^n(1-(-z)^j)^{C_j^{(n)}},
 \qquad k_n/n\to c\in(0,1).
 \label{eq:intermediate-hook}
\end{equation}
Each cycle contributes its determinant factor, and removing the
invariant line divides the product by $1+z$. Unlike the fixed-tail
formula, this coefficient depends on more than a fixed finite
cycle vector. To study its subsequential laws, one can try to
control the growing-degree coefficient as a function of the whole
cycle vector, including cycles whose lengths grow with $n$.
Once an ordinary character limit is obtained, no additional
square-character estimate is needed to transfer it to the
randomized-basis statistic.

For two rows, the subset-module Pieri decomposition gives the
character as the difference of the coefficients of $z^{k_n}$ and
$z^{k_n-1}$ in $\prod_j(1+z^j)^{C_j^{(n)}}$. These formulas give
starting points for character laws beyond fixed tails and
sublinear largest hooks. Thus the remaining diagram questions
concern the ordinary character that accompanies the Gaussian
fluctuation; the criterion for the centered Gaussian part has
already been established.

\newpage

\section*{Acknowledgements}
The author thanks Steve Evans for proposing this problem, namely
by pointing him to Tsou~\cite{TsouRandomized} to understand the 
problem statement, as well as to Collins, Matsumoto, and Novak
~\cite{CollinsMatsumotoNovak} for an introduction to 
Weingarten calculus. The author is grateful to Evans for the 
helpful comments and suggestions, such as to  relate eventual 
cycle-polynomial character formulas to  stabilization of Young
diagrams after deletion of the first row or first column, and 
for proposing the use of the Murnaghan–-Nakayama rule in the 
converse argument of Proposition~\ref{prop:exact-tail-equivalence}.

\end{document}